\documentclass[10pt]{article}

\usepackage{amssymb} \usepackage{amsmath}
\usepackage{latexsym} \usepackage{graphicx}

\newtheorem{Thm}{Theorem}
\newtheorem{Cor}{Corollary}
\newtheorem{Lem}{Lemma}

  \def\bpf{\noindent{\sc
Proof.~}} \def\epf{\quad \quad \quad \qed} \def\qed{\vbox{\hrule
\hbox{\vrule\hbox to 5pt{\vbox to 8pt{\vfil}\hfil}\vrule}\hrule}}
 \def\({\left(} \def\){\right)}

\DeclareRobustCommand{\binom}{\frbinom{}}
\def\frbinom#1#2#3{{#1{#2\atopwithdelims()#3}}}

\title{On the Likelihood of Comparability in Bruhat Order}
\author{{\sc Adam Hammett} and {\sc Boris Pittel}\\ {\small Department
of Mathematics} \\ {\small The Ohio State University} \\ {\small 231
W. 18$^\mathrm{th}$ Ave., Columbus, OH 43210} \\ {\small \texttt{
hammett@math.ohio-state.edu$^{*}$}} \\ {\small \texttt{
bgp@math.ohio-state.edu$^{**}$}}} \date{}

\begin{document}

\maketitle

\def\thefootnote{}
\footnote{ \\
\noindent $^{*}$Supported in part by NSF grant DMS-0104104.

\noindent $^{**}$Supported in part by NSF grants DMS-0104104 and DMS-0406024. }

\begin{abstract}

\noindent Two permutations of $\left[ n \right]$ are comparable in the
Bruhat order if one is closer, in a natural way, to the identity
permutation, $1 \, 2 \, \cdots \, n$, than the other. We show that the
number of comparable pairs is of order $\(n!\)^2/n^2$ at most, and
$\(n!\)^2\(0.708\)^n$ at least. For the \emph{weak} Bruhat order, the
corresponding bounds are $\(n!\)^2\(0.362\)^n$ and
$\(n!\)^2\prod_{i=1}^n \(H\(i\)/i\)$, where $H\(i\):=\sum_{j=1}^i
1/j$. In light of numerical experiments, we conjecture that for each
order the upper bound is qualitatively close to the actual number of
comparable pairs.
\end{abstract}

\bigskip
\thispagestyle{empty}

\begin{center}
{\Large {\bf Introduction}}
\end{center}

Let $n \geq 1$ be an integer. Two permutations of
$[n]:=\left\{1,\dots,n\right\}$ are comparable in the \emph{Bruhat
order} if one can be obtained from the other by a sequence of
transpositions of pairs of elements forming an inversion. Here is
a precise definition of the (strong)
Bruhat order on the set of permutations $S_n$
(see Stanley \cite[p. 172, ex. 75. a.]{S}, Humphreys
\cite[p. 119]{H}). If $\omega = \omega\( 1
\) \cdots \omega \( n \)\in S_n$, then a \textit{reduction}
of $\omega$ is a permutation obtained from $\omega$ by interchanging
some $\omega \( i \)$ with some $\omega \( j \)$ provided $i<j$ and
$\omega \( i \)>\omega\( j \)$. We say that $\pi \leq \sigma$ in the
Bruhat order if there is a chain $\sigma=\omega_1\to\omega_2\to\cdots \to \omega_s=
\pi$, where each $\omega_t$ is a reduction of $\omega_{t-1}$. The number
of inversions in $\omega_t$ strictly decreases with $t$. Indeed, one can show that if 
$\omega_2$ is a reduction of $\omega_1$ via the interchange 
$\omega_1(i)\leftrightarrow \omega_1(j)$, $i<j$, then

\begin{eqnarray*}
& \text{inv}(\omega_1)=\text{inv}(\omega_2)+2N(\omega_1)+1,\\
&N(\omega_1):=|\{ k \, : \, i<k<j,\, \omega_1(i)>\omega_1(k)>\omega_1(j)\}|;
\end{eqnarray*}

\noindent here $\text{inv}(\omega_1)$, say, is the number of inversions in $\omega_1$ 
(see Bj\"orner and Brenti \cite{BB1}). Figure 1 below illustrates this poset on $S_3$ and $S_4$.

The definition of the Bruhat order is very transparent, and yet
deciding for given $\pi,\sigma$ whether $\pi\le\sigma$ is
computationally difficult, even for smallish $n$. Fortunately, the historically
first definition of Bruhat comparability, due to Ehresmann (1934) \cite{E}, 
is actually an efficient algorithm to check if $\pi\le\sigma$ 
(see also Knuth \cite{K}, \cite{BB1}). The Ehresmann ``tableau criterion''
states that $\pi \le \sigma $ if and only if
$\pi_{i,j}\le \sigma_{i,j}$ for all $1\le i\le j\le n-1$, where $\pi
_{i,j}$ and $\sigma_{i,j}$ are the $i$-th entry in the increasing rearrangement of $\pi \(
1 \) ,\dots ,\pi \( j \) $ and of $\sigma\(1\),\dots,\sigma\(j\)$. These arrangements
form two staircase tableaux, hence the term ``tableau criterion''.
For example, $41523>21534$ is verified by element-wise
comparisons of the two tableaux

\[
\begin{array}{c}
\begin{array}{cccc}
1 & 2 & 4 & 5 \\ 1 & 4 & 5 & \\ 1 & 4 & & \\ 4 & & &
\end{array}
\qquad \qquad
\begin{array}{cccc}
1 & 2 & 3 & 5 \\ 1 & 2 & 5 & \\ 1 & 2 & & \\ 2 & & &
\end{array}
\end{array}.
\]

\noindent Then, Deodhar \cite{D} extended the Bruhat order notion to
other Coxeter groups. It is well-known that Ehresmann's criterion is
equivalent to the $\(0,1\)$-matrix criterion. It involves comparing
the number of $1$'s contained in certain submatrices of the
$\(0,1\)$-permutation matrices representing $\pi$ and $\sigma$ (see
B\'ona \cite{Bo}, \cite{BB1}). Later, Bj\"orner and Brenti \cite{BB2}
were able to improve on the result of \cite{E}, giving a tableau
criterion that requires fewer operations. Very recently, Drake,
Gerrish and Skandera \cite{DGS} have found two new comparability
criteria, involving totally nonnegative polynomials and the Schur
functions respectively. We are aware of other criteria (see Bj\"orner
\cite{Bj}, Fulton \cite[pp. 173-177]{Fu}, Lascoux and Sch\"utzenberger
\cite{LS}, \cite{D}), but we found the $\(0,1\)$-matrix and Ehresmann
criteria most amenable to probabilistic study.

\bigskip
\begin{center}
\includegraphics[scale=0.7]{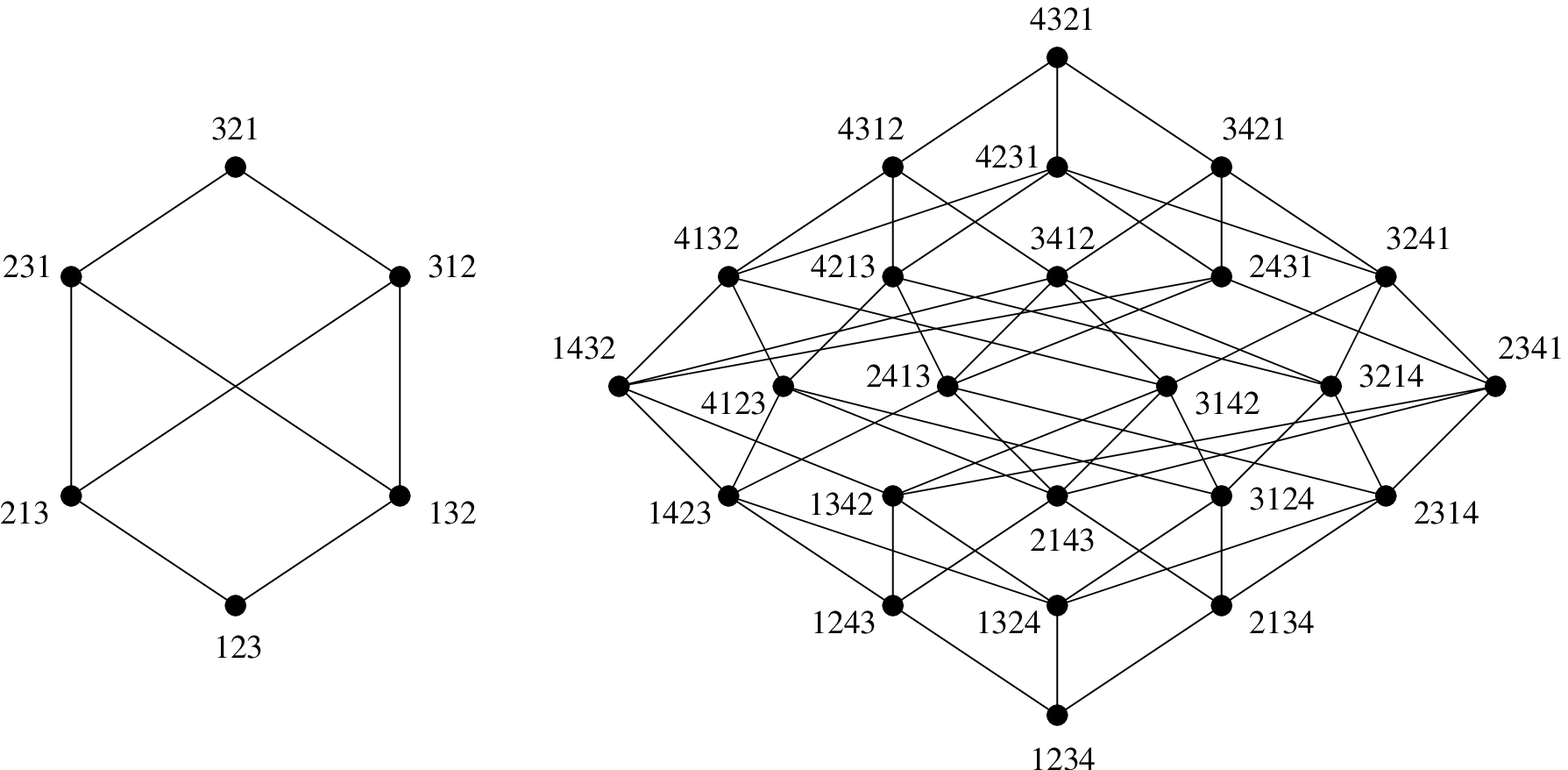} \\ {\small {\bf Figure 1:}
The Bruhat order on $S_3$ and $S_4$.}
\end{center}
\bigskip

The $\(0,1\)$-matrix criterion for Bruhat order on $S_n$ says that for
$\pi ,\sigma \in S_n$, $\pi \le \sigma$ if and only if for all $i,j
\le n$, the number of $\pi \(1\),\dots,\pi \(i \)$ that are at most
$j$ exceeds (or equals) the number of $\sigma \(1\),\dots,\sigma
\(i\)$ that are at most $j$ (see \cite{Bo} for this version). It is
referred to as the $\(0,1\)$-matrix criterion because of the following
recasting of this condition: let $M\(\pi\)$, $M(\sigma)$ be the
permutation matrices corresponding to $\pi$, $\sigma$, so that for
instance the $\(i,j\)$-entry of $M(\pi)$ is $1$ if $\pi\(j\)=i$ and
$0$ otherwise. Here, we are labeling columns $1,2,\dots,n$ when
reading from \emph{left to right}, and rows are labeled $1,2,\dots,n$ when
reading from \emph{bottom to top} so that this interpretation is like placing
ones at points $\(i,\pi\(i\)\)$ of the $n \times n$ integer lattice and 
zeroes elsewhere. Denoting submatrices of
$M\(\cdot\)$ corresponding to rows $I$ and columns $J$ by
$M\(\cdot\)_{I,J}$, this criterion says that $\pi \le \sigma$ if and
only if for all $i,j \le n$, the number of ones in
$M\(\pi\)_{[i],[j]}$ is at least the number of ones in
$M\(\sigma\)_{[i],[j]}$ (see \cite{DGS} for this version).

An effective way of visualizing this criterion is to imagine the
matrices $M\(\pi\)$ and $M\(\sigma\)$ as being superimposed on one
another into a single matrix, $M\(\pi,\sigma\)$, with the ones for
$M\(\pi\)$ represented by $\times$'s (``crosses''), the ones for
$M\(\sigma\)$ by $\circ$'s (``balls'') and the zeroes for both by
empty entries. Note that some entries of $M\(\pi,\sigma\)$ may be
occupied by both a cross and a ball. Then the $\(0,1\)$-matrix
criterion says that $\pi \le \sigma$ if and only if every southwest
submatrix of $M\(\pi,\sigma\)$ contains at least as many crosses as
balls. Here, in the notation above, a \emph{southwest submatrix} is a
submatrix $M\(\pi,\sigma\)_{[i],[j]}$ of $M\(\pi,\sigma\)$ for some
$i,j \le n$. It is clear that we could also check $\pi \le \sigma$ by
checking that crosses are at least as numerous as balls in every
northeast submatrix of $M\(\pi,\sigma\)$. Likewise, $\pi \le \sigma$
if and only if balls are at least as numerous as crosses in every
northwest submatrix of $M\(\pi,\sigma\)$, or similarly balls are at
least as numerous as crosses in every southeast submatrix of
$M\(\pi,\sigma\)$. Parts of all four of these equivalent conditions
will be used in our proofs. As a quick example, with $\pi=21534$ and
$\sigma=41523$, $\pi < \sigma$ is checked by examining southwest
submatrices of $M\(\pi,\sigma\)$ in Figure 2 below. Also, the
superimposing of $M\(\pi\)$ with $M\(\sigma\)$ to form
$M\(\pi,\sigma\)$ is illustrated in this figure.

\bigskip
\begin{center}
\includegraphics[scale=0.7]{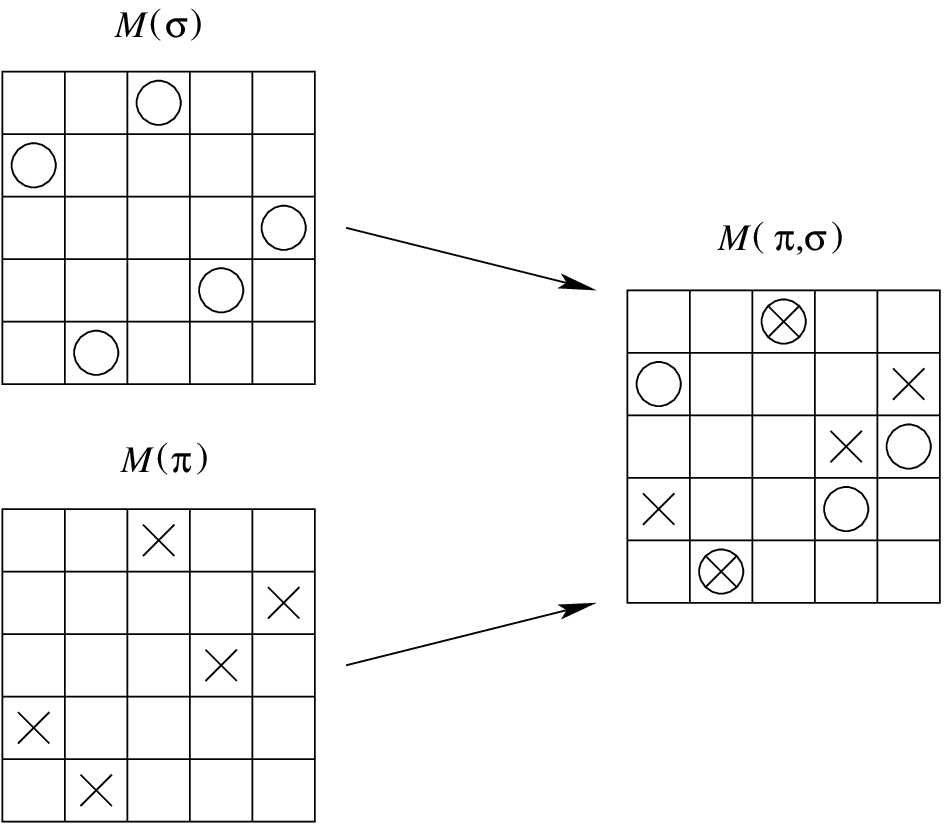} \\ {\small {\bf Figure 2:}
Superimposing $M(\pi)$ and $M(\sigma)$ to form $M(\pi,\sigma)$.}
\end{center}
\bigskip

\noindent In this note, we use the $\(0,1\)$-matrix and the Ehresmann criteria
to obtain upper and lower bounds for the number of pairs $(\pi,\sigma)$ with 
$\pi\le\sigma$.

\bigskip

\begin{Thm} \label{ThmSBr} Let $n \ge 1$ be an integer, and let $\pi,\sigma \in S_n$
be selected independently and uniformly at random. Then there exist
universal constants $c_1,c_2>0$ such that

\[
c_1\(0.708\)^n \le P\( \pi \leq \sigma \) \le c_2/n^2.
\]

\end{Thm}

\bigskip

\noindent Equivalently, the number of pairs $(\pi,\sigma)$ with $\pi \le \sigma$
is sandwiched between $c_1(0.708)^n(n!)^2$ and $c_2n^{-2}(n!)^2$. The
lower bound follows from a sufficient condition derived from the
$(0,1)$-matrix criterion, and a computer-aided tabulation of an
attendant function of a smallish integer argument.
Empirical estimates based on generating pairs of random
permutations suggest that $P\(\pi \le \sigma\)$ is of order $n^{-\(2+\delta\)}$, 
for $\delta$ close to $0.5$ from above. So apparently it is the upper bound which
comes close to the true proportion $P\( \pi \le \sigma\)$. It is certain
that the constant $0.708$ can be further improved, but we do not know if 
our method could be extended to deliver a lower bound $(1-o(1))^n$.
A lower bound $n^{-a}$, a qualitative match of the
upper bound, seems out of sight presently.

Then we turn to the modified order on $S_n$, the \textit{weak} Bruhat
order ``$\preceq$''. Here $\pi\preceq\sigma$ if there is a chain
$\sigma=\omega_1\to\omega_2\to\cdots\to\omega_s=\pi$, where
each $\omega_t$ is a \textit{simple reduction} of $\omega_{t-1}$,
i.e. obtained from $\omega_{t-1}$ by transposing two \emph{adjacent} elements
$\omega_{t-1}(i)$, $\omega_{t-1}(i+1)$ with $\omega_{t-1}(i)>\omega_{t-1}(i+1)$. 
Since at each step the number of inversions decreases by $1$, all chains
connecting $\sigma$ and $\pi$ have the same length.

\bigskip

\begin{Thm} \label{ThmWBr} Let $P_n^*:=P\(\pi\preceq
\sigma\)$.  Then $P_n^*$ is submultiplicative, i.e. $P_{n_1+n_2}^*\le
P_{n_1}^*P_{n_2}^*$. Consequently there exists $\rho=\lim
\sqrt[n]{P_n^*}$. Furthermore, there exists an absolute constant $c>0$
such that

\[
\prod_{i=1}^n \( H\(i\)/i\) \le P_n^* \le c\(0.362\)^n,
\]

\noindent where $H\(i\) := \sum_{j=1}^i 1/j$. Consequently, $\rho \le
0.362$.

\end{Thm}

\bigskip

\noindent The proof of the upper bound is parallel to that of Theorem
\ref{ThmSBr}, lower bound, while the lower bound follows from an
inversion set criterion inspired by discussion of the weak Bruhat
order by Berge \cite{Be}.  Empirical estimates indicate that
$\rho$ is close to $0.3$.  So here too the upper bound
seems to be qualitatively close to the actual probability $P_n^*$.
And our lower bound, though superior to the trivial bound $1/n!$, is
decreasing superexponentially fast with $n$, which makes us believe
that there ought to be a way to vastly improve it. Paradoxically, it
is the lower bound that required a deeper combinatorial insight. We
show first that the number of $\pi$'s below (or equal to) $\sigma$
equals $e({\cal P})$, the total number of linear extensions of ${\cal
P}={\cal P}(\sigma)$, the poset induced by $\sigma$.  (The important
notion of ${\cal P}(\sigma)$ was brought to our attention by Sergey
Fomin \cite{Fo}.) And then we prove that for each poset ${\cal P}$ of
cardinality $n$,

\[
e\({\cal P}\) \ge n! \Big/ \prod_{i\in {\cal P}} d\(i\),
\]

\noindent where $d\(i\) := | \left\{ j \in {\cal P} \, : \, j \le i
\text{ in } {\cal P} \right\} |$. The final step is based on the
independence of sequential ranks in the uniform permutation. Reducing
the gap between the bounds will probably require a better
understanding of ${\cal P}(\sigma)$ for a typical permutation
$\sigma$.

In conclusion we mention two papers, \cite{P1} and \cite{P2}, where
the ``probability-of-comparability'' problems were solved for the poset of 
integer partitions of $n$ under dominance order, and
for the poset of set partitions of $\left[n\right]$ ordered by refinement.

\bigskip

\begin{center}
{\Large {\bf Proof of Theorem \ref{ThmSBr}, upper bound.}}
\end{center}

We need to show that

\[
P\( \pi \le \sigma \) = O \( n^{-2} \).
\]

\noindent The argument is based on the $(0,1)$-matrix criterion.  We
assume that $n$ is even. Only minor modifications are necessary for
$n$ odd.

\bigskip

STEP 1. The $(0,1)$-matrix criterion requires that a set of $n^2$
conditions are met.  The challenge is to select a subset of those
conditions which meets two conflicting demands. It has to be
sufficiently simple so that we can compute (estimate) the probability
that the random pair $(\pi,\sigma)$ satisfies all the chosen
conditions. On the other hand, collectively these conditions need to
be quite stringent for this probability to be $o(1)$. In our first
advance we were able (via Ehresmann's criterion) to get a bound $O(n^{-1/2})$ by 
using about $2n^{1/2}$ conditions.  We are about to describe a set of $2n$
conditions that does the job.

Let us split the matrices $M\(\pi,\sigma\)$, $M\(\pi\)$ and
$M\(\sigma\)$ into 4 submatrices of equal size $n/2 \times n/2$ -- the
southwest, northeast, northwest and southeast corners, denoting them
$M_{sw}\(\cdot\)$, $M_{ne}\(\cdot\)$, $M_{nw}\(\cdot\)$ and
$M_{se}\(\cdot\)$ respectively. In the southwest corner
$M_{sw}\(\pi,\sigma\)$, we restrict our attention to southwest
submatrices of the form $i \times n/2$, $i=1,\dots,n/2$. If $\pi \le
\sigma$, then as we read off rows of $M_{sw}\(\pi,\sigma\)$ from
bottom to top keeping track of the total number of balls and crosses
encountered thus far, at any intermediate point we must have at least
as many crosses as balls. Let us denote the set of pairs
$\(\pi,\sigma\)$ such that this occurs by $\mathcal{E}_{sw}$. We draw
analogous conclusions for the northeast corner, reading rows from top
to bottom, and we denote by $\mathcal{E}_{ne}$ the set of pairs
$\(\pi,\sigma\)$ satisfying this condition.

Similarly, we can read columns from left to right in the northwest
corner, and here we must always have at least as many balls as
crosses.  Denote the set of these pairs $\(\pi,\sigma\)$ by
$\mathcal{E}_{nw}$. The same condition holds for the southeast corner
when we read columns from right to left. Denote the set of these pairs
$\(\pi,\sigma\)$ by $\mathcal{E}_{se}$. Letting $\mathcal{E}$ denote
the set of pairs $\(\pi,\sigma\)$ satisfying all four of the
conditions above, we get

\begin{equation*}
\left\{ \pi \le \sigma \right\} \subseteq \mathcal{E}=\mathcal{E}_{sw}
\cap \mathcal{E}_{ne} \cap \mathcal{E}_{nw} \cap \mathcal{E}_{se}.
\end{equation*}

\noindent Pairs of permutations in $\mathcal{E}$ satisfy $2n$ of the
$n^2$ conditions required by the $\(0,1\)$-matrix criterion. And
unlike the set $\{\pi\le\sigma\}$, we are able to compute
$|\mathcal{E}|$, and to show that
$P(\mathcal{E})=(n!)^{-2}|\mathcal{E}|=O(n^{-2})$.  Figure 3 below is
a graphical visualization of the reading-off process that generates
the restrictions defining the set $\mathcal{E}$.

\bigskip
\begin{center}
\includegraphics[scale=0.7]{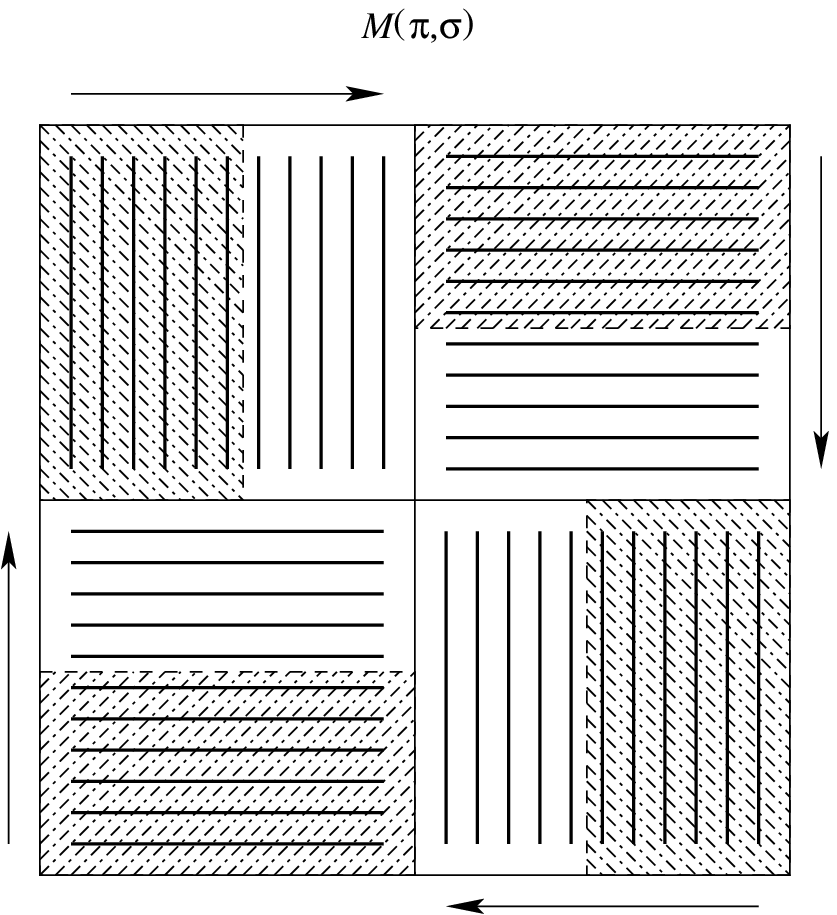} \\ {\small {\bf Figure 3:}
Finding a necessary condition for $\pi \le \sigma$.}
\end{center}
\bigskip

If a row (column) of a submatrix $M(\pi)_{I,J}$ ($M(\sigma)_{I,J}$ resp.) 
contains a marked entry, we say that it {\it supports\/} the submatrix. Clearly 
the number of supporting rows (columns) equals the number of marked
entries in $M(\pi)_{I,J}$ ($M(\sigma)_{I,J}$ resp.).  Now, given $\pi$, $\sigma$, 
let $M_1=M_1(\pi)$, $M_2=M_2(\sigma)$ denote the total number of rows that
support $M_{sw}(\pi)$ and $M_{sw} (\sigma)$ respectively. Then
$M_{nw}(\pi)$, $M_{nw}(\sigma)$ are supported by $M_3=n/2-M_1$ columns
and by $M_4=n/2-M_2$ columns respectively. The same holds for the
southeastern corners of $M(\pi)$ and $M(\sigma)$. Obviously the
northeastern submatrices of $M(\pi)$ and $M(\sigma)$ are supported by
$M_1$ rows and $M_2$ rows respectively. Then we have

\begin{eqnarray} \label{ballotEqu}
&P\( \mathcal{E}\)=\sum\limits_{m_1,m_2}P\(\mathcal{E}\cap
{\cal A}\(m_1,m_2\)\),\\
&{\cal A}\(m_1,m_2\):=\left\{(\pi,\sigma):
\,M_1=m_1,M_2=m_2 \right\}. \nonumber
\end{eqnarray}
 
\noindent Clearly $\mathcal{E}\cap
{\cal A}\(m_1,m_2\)=\emptyset$ if $m_1< m_2$. We
claim that, for $m_1\ge m_2$,

\begin{equation} \tag{\dag}
P\(\mathcal{E}\cap {\cal A}\(m_1,m_2\)\)= 
\left[ \frac{(m_1-m_2+1)(n/2+1)}{(m_3+1)(m_1+1)} \right]^4 \cdot
\frac{\prod_{i=1}^4 \binom{n/2}{m_i}}{\binom{n}{n/2}^2}.
\end{equation}

\noindent Here and below $m_3:=n/2-m_1$ and $m_4:=n/2-m_2$ stand for generic
values of $M_3$ and $M_4$ in the event ${\cal A}\(m_1,m_2\)$.

To prove (\dag), let us count the number of pairs $(\pi,\sigma)$ in
$\mathcal{E}\cap {\cal A}\(m_1,m_2\)$. First consider the southwest corner, 
$M_{sw}(\pi,\sigma)$. Introduce $L_1=L_1\(\pi,\sigma\)$, the number of rows 
supporting both $M_{sw}(\pi)$ and $M_{sw}(\sigma)$. So $L_1$ is the number 
of rows in the southwest corner $M_{sw}(\pi,\sigma)$ containing both a cross 
and a ball. Suppose that we are on the event $\{L_1=\ell_1\}$. 
We choose $\ell_1$ rows to support both $M_{sw}\(\pi\)$ and $M_{sw}\(\sigma\)$
from the $n/2$ first rows. Then, we choose $\(m_1-\ell_1+m_2-\ell_1\)$
more rows from the remaining $(n/2-\ell_1)$ rows.  Each of these
secondary rows is to support either $M_{sw}(\pi)$ or $M_{sw}(\sigma)$,
but not both. This step can be done in

\[
\binom{n/2}{\ell_1}\binom{n/2-\ell_1}{m_1-\ell_1+m_2-\ell_1}
\]

\noindent ways. Next, we partition the set of
$\(m_1-\ell_1+m_2-\ell_1\)$ secondary rows into two row subsets of
cardinality $(m_1-\ell_1)$ and $(m_2-\ell_1)$ that will support
$M_{sw}\(\pi\)$ and $M_{sw}(\sigma)$, accompanying the $\ell_1$
primary rows supporting both submatrices. We can visualize each of the
resulting row selections as a subsequence of $(1,\dots,n/2)$ which is
a disjoint union of two subsequences, with $\ell_1$ and $(
m_1-\ell_1+m_2-\ell_1 )$ elements respectively, with each element of
the subsequence marked by a cross and/or ball dependent on which
submatrix is supported by this element, i.e. the row represented by
this element.  The condition $\mathcal{E}_{sw}$ is equivalent to the
restriction: moving along the subsequence from left to right, at each
point the number of crosses is not to fall below the number of balls.
Obviously, no double-marked element can cause violation of this
condition.  Thus, our task is reduced to determination of the number
of $(m_1-\ell_1 +m_2 -\ell_1)$-long sequences of $m_1-\ell_1$ crosses
and $m_2-\ell_1$ balls such that at no point the number of crosses is
strictly less than the number of balls.  By the classic ballot theorem
(see Takacs \cite[pp. 2-7]{T}), the total number of such sequences
equals

\[
\frac{\(m_1-\ell_1+1\)-\(m_2-\ell_1\)}{\(m_1-\ell_1+1\)+\(m_2-\ell_1\)}
\binom{m_1-\ell_1+m_2-\ell_1+1}{m_1-\ell_1+1}=\frac{m_1-m_2+1}{m_1-\ell_1+1}
\binom{m_1-\ell_1+m_2-\ell_1}{m_1-\ell_1}.
\]

\noindent The second binomial coefficient is the total number of
$(m_1-\ell_1+m_2-\ell_1)$-long sequences of $(m_1-\ell_1)$ crosses and
$(m_2- \ell_1)$ balls. So the second fraction is the probability that
the sequence chosen uniformly at random among all such sequences meets
the ballot theorem condition. The total number of ways to designate
the rows supporting $M_{sw}(\pi)$ and $M_{sw}(\sigma)$, subject to the
condition $\mathcal{E}_{sw}$, is the product of two counts, namely

\begin{equation*} \begin{split}
\binom{n/2}{\ell_1}\binom{n/2-\ell_1}{m_1-\ell_1+m_2-\ell_1}
&\binom{m_1-\ell_1+m_2-\ell_1}{m_1-\ell_1}\frac{m_1-m_2+1}{m_1-\ell_1+1} \\
&=\frac{m_1-m_2+1}{n/2-m_2+1}\binom{n/2}{m_2}\binom{m_2}{\ell_1}
\binom{n/2-m_2+1}{m_1-\ell_1+1}.
\end{split} \end{equation*}

\noindent Summing this last expression over all $\ell_1 \le m_2$, we obtain

\begin{equation} \label{ballotEquFactor} \begin{split}
\frac{m_1-m_2+1}{n/2-m_2+1}\binom{n/2}{m_2}\sum_{\ell_1 \le m_2}\binom{m_2}{\ell_1}
\binom{n/2-m_2+1}{m_1-\ell_1+1}
&= \frac{m_1-m_2+1}{n/2-m_2+1}\binom{n/2}{m_2}\binom{n/2+1}{m_1+1} \\
&= \frac{(m_1-m_2+1)(n/2+1)}{(n/2-m_2+1)(m_1+1)}\binom{n/2}{m_1}\binom{n/2}{m_2}.
\end{split} \end{equation}

\noindent Here, in the first equality, we have used the binomial theorem. The product of the 
two binomial coefficients in the final count (\ref{ballotEquFactor}) is the total number of row 
selections from the first $n/2$ rows, $m_1$ to contain crosses and $m_2$ to contain balls. 
So the fraction preceding these two binomial factors is the probability that a particular row 
selection chosen uniformly at random from all such row selections satisfies our ballot 
condition ``crosses never fall below balls''. Equivalently, by the very derivation, the 
expression (\ref{ballotEquFactor})
is the total number of paths $(X(t),Y(t))_{0 \le t \le n/2}$ on the square lattice connecting $(0,0)$ 
and $(m_1,m_2)$ such that $X(t+1)-X(t),\, Y(t+1)-Y(t) \in \{0,1\}$, and $X(t)\ge Y(t)$ for every $t$. 
(To be sure, if $X(t+1)-X(t)=1$ and $Y(t+1)-Y(t)=1$, the corresponding move is a combination of  horizontal and vertical unit moves.)

Likewise, we consider the northeast corner, $M_{ne}(\pi,\sigma)$. We introduce $L_2=L_2(\pi,\sigma)$, 
the number of rows in $M_{ne}(\pi,\sigma)$ containing both a cross and a ball. By initially 
restricting to the event $\{L_2=\ell_2\}$, then later summing over all $\ell_2 \le m_2$, 
we obtain another factor (\ref{ballotEquFactor}). Analogously, a third and fourth factor 
(\ref{ballotEquFactor}) comes from considering {\it columns\/} 
in the northwest and southeast corners, $M_{nw}(\pi,\sigma)$ and $M_{se}(\pi,\sigma)$. 
Importantly, the row selections for the southwest and the northeast
submatrices do not interfere with the column selections for the
northwest and the southeast corners. So by multiplying these four factors (\ref{ballotEquFactor}) 
we obtain the total number of row and column selections on the event ${\cal A}(m_1,m_2)$ 
subject to all four restrictions defining ${\cal E}$!

Once such a row-column selection has been made, we have determined
which rows and columns support the four submatrices of $M(\pi)$ and
$M(\sigma)$.  Consider, for instance, the southwest corner of
$M(\pi)$. We have selected $m_1$ rows (from the first $n/2$ rows)
supporting $M_{sw}(\pi)$, and we have selected $m_3$ columns (from the
first $n/2$ columns) supporting $M_{nw}(\pi)$. Then it is the
remaining $n/2-m_3=m_1$ columns that support $M_{sw}(\pi)$. The number
of ways to match these $m_1$ rows and $m_1$ columns, thus to determine
$M_{sw} (\pi)$ completely, is $m_1!$. The northeast corner contributes
another $m_1!$, while each of the two other corners contributes
$m_3!$, whence the overall matching factor is $(m_1!m_3!)^2$. The
matching factor for $\sigma$ is $(m_2!m_4!)^2$. Multiplying the number
of admissible row-column selections by the resulting
$\prod_{i=1}^4(m_i!)^2$ and dividing by $(n!)^2$, we obtain

\[
P\({\cal E}\cap {\cal A}\(m_1,m_2\)\)
= \left[\frac{(m_1-m_2+1)(n/2+1)}{(n/2-m_2+1)(m_1+1)}\binom{n/2}{m_1}\binom{n/2}{m_2}\right]^4 \cdot
\frac{\prod_{i=1}^4(m_i!)^2}{(n!)^2},
\]

\noindent which is equivalent to (\dag). Figure 4 below is a graphical explanation of 
this matching factor. In it, we show the matrix $M(\pi)$ in a case when in the southwest and the
northeast squares $\pi$ is supported by the bottom $m ( =m_1 )$ and the top $m$
rows respectively; likewise, in the northwest and the southeast squares $\pi$ is
supported by the $n/2-m$ leftmost and the $n/2-m$ rightmost columns respectively.

\bigskip
\begin{center}
\includegraphics[scale=0.9]{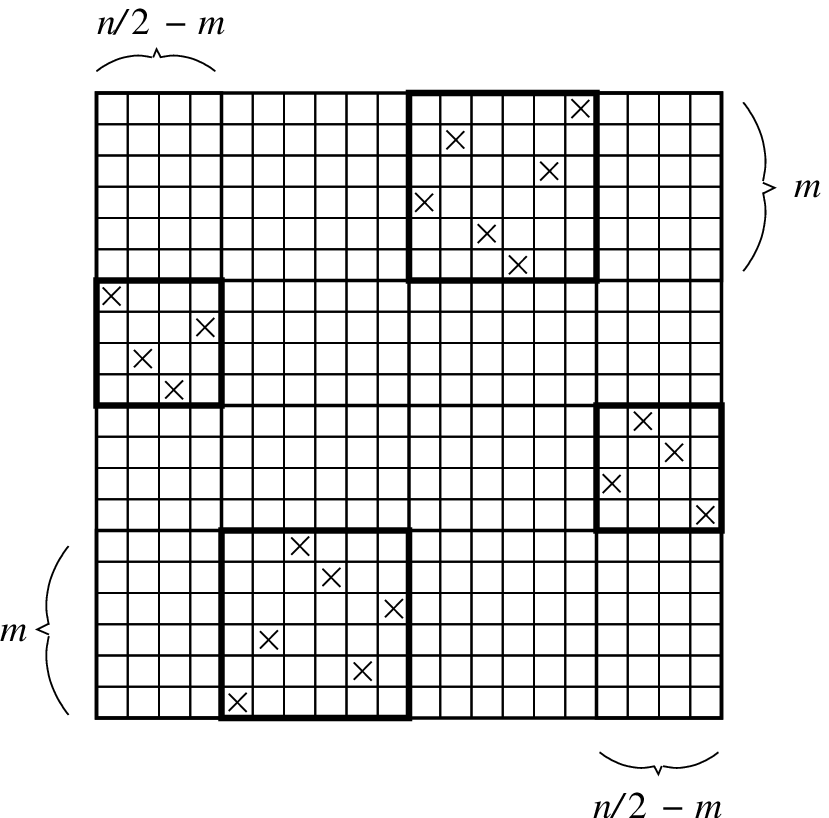} \\ 
{\small {\bf Figure 4:}
Selection of first $m=m_1$ ($n/2-m$ resp.) rows (columns resp.) in corners to support $M(\pi)$.}
\end{center}
\bigskip

STEP 2. Let us show that (\ref{ballotEqu}) and (\dag) imply

\begin{equation*} \tag{\ddag} \begin{split}
P\(\mathcal{E}\)\le E\left[
\frac{\(M_1-M_2+1\)^4\(n/2+1\)^4}{\(n/2-M_2+1\)^4\(M_1+1\)^4}
\right].
\end{split} \end{equation*}

\noindent First, $M_1$ and $M_2$ are independent with

\begin{equation*}
P\(M_i=m_i\)=\frac{\binom{n/2}{m_i}^2}{\binom{n}{n/2}}, \quad i=1,2.
\end{equation*}

\noindent Indeed, $M_i$ obviously
equals the cardinality of the intersection with $[n/2]$ of a uniformly random 
subset of size $n/2$ from $[n]$, which directly implies these formulas. Thus, each $M_i$
has the hypergeometric distribution with parameters $n/2,n/2,n/2$; in
other words, $M_i$ has the same distribution as the number of red
balls in a uniformly random sample of $n/2$ balls from an urn
containing $n/2$ red balls and $n/2$ white balls. By the independence of $M_1$ and $M_2$, 
we obtain

$$
P(M_1=m_1,M_2=m_2)=\frac{\binom{n/2}{m_1}^2\binom{n/2}{m_2}^2}{\binom{n}{n/2}^
2}.
$$

\noindent It remains to observe that (\ref{ballotEqu}) and (\dag) imply

\begin{equation*} \begin{split}
P({\cal E})&=\sum_{m_1 \ge m_2} \frac{\(m_1-m_2+1\)^4\(n/2+1\)^4}{\(n/2-m_2+1\)^4\(m_1+1\)^4} 
\cdot P(M_1=m_1,M_2=m_2) \\
&\le \sum_{m_1,m_2} \frac{\(m_1-m_2+1\)^4\(n/2+1\)^4}{\(n/2-m_2+1\)^4\(m_1+1\)^4} 
\cdot P(M_1=m_1,M_2=m_2) \\
&=E\left[ \frac{\(M_1-M_2+1\)^4\(n/2+1\)^4}{\(n/2-M_2+1\)^4\(M_1+1\)^4} \right],
\end{split} \end{equation*}

\noindent and (\ddag) is proved.

\bigskip

STEP 3. The advantage of (\ddag) is that it allows us to use
probabilistic tools exploiting the independence of the random
variables $M_1$ and $M_2$. Typically the $M_i$'s are close to $n/4$, while 
$|M_1-M_2|$ is of order $n^{1/2}$ at most. So, in view
of (\ddag) we expect that $P(\mathcal{E}) =O(n^{-2})$.

We now make this argument rigorous. First of all, by the
``sample-from-urn'' interpretation of $M_i$,

\begin{equation}\label{E[M]}
E\left[ M_i \right]
=\frac{n}{2}\frac{\binom{n-1}{n/2-1}}{\binom{n}{n/2}}= n/4.
\end{equation}

\noindent Then (see Janson et al. \cite{JLR}) the probability
generating function of $M_i$ is dominated by that of
$\text{Bin}(n,1/4)$, and consequently for each $t \ge 0$ we have

\[
P \( |M_i-n/4| \ge t \) = O\( \exp \(-4t^2/n \)\).
\]

\noindent Hence, setting $t=n^{2/3}$ we see that

\begin{equation*}
P\( n/4-n^{2/3} < M_i < n/4+n^{2/3} \)\ge 1-e^{-cn^{1/3}},
\end{equation*}

\noindent for some absolute constant $c>0$. Introduce the
event

$$
A_n=\bigcap_{i=1}^2\left\{n/4-n^{2/3}< M_i< n/4+n^{2/3}\right\}.
$$

\noindent Combining the estimates for $M_i$, we see that for some 
absolute constant $c_1>0$,

$$ 
P(A_n)\ge 1-e^{-c_1n^{1/3}}.
$$

\noindent Now the random variable in (\ddag), call it $X_n$, is
bounded by $1$, and on the event $A_n$, within a factor of
$1+O(n^{-1/3})$,

$$ 
X_n=\left(\frac{4}{n}\right)^8(M_1-M_2+1)^4\(n/2+1\)^4.
$$

\noindent Therefore

\begin{equation*}\begin{split}
P\(\mathcal{E}\) &\le \left(\frac{5}{n}\right)^8\(n/2+1\)^4 E
\left[(M_1-M_2+1)^4\right] +O\left(e^{-c_1n^{1/3}}\right).
\end{split} \end{equation*} 

\noindent It remains to prove that this expected value is
$O(n^2)$. Introduce $\overline{M}_i=M_i-E[M_i]$, $i=1,2$. Then

\begin{equation*}
(M_1-M_2+1)^4=(\overline{M}_1-\overline{M}_2+1)^4\le
27(\overline{M}_1^4+\overline {M}_2^4+1),
\end{equation*}

\noindent as

\[
(a+b+c)^2\le 3(a^2+b^2+c^2).
\]

\noindent We now demonstrate that $E[\overline{M}_i^4]=O(n^2)$. 
To reduce computations to the minimum we use the fact that, as a
special instance of the hypergeometrically distributed random
variable, $M_i$ has the same distribution as the sum of $n/2$
independent Bernoulli variables $Y_j\in\{0,1\}$ (see Vatutin and
Mikhailov \cite{VM}, alternatively \cite[p. 30]{JLR}). So denoting
$p_j=P(Y_j=1)$, $q_j=P(Y_j=0)$, we get a product formula for the
moment generating function of $\overline{M}_i$:

\begin{equation*}\begin{split}
E\left[e^{u\overline{M}_i}\right]&=\prod_{j=1}^{n/2}E\left[e^{u(Y_j-E[Y_j])}\right]
=\prod_{j=1}^{n/2}\left(p_je^{uq_j}+q_je^{-up_j}\right)\\
&=\prod_{j=1}^{n/2}\left(1+\frac{u^2}{2}p_jq_j+\frac{u^3}{3!}p_jq_j(q_j-p_j)
+\frac{u^4}{4!}(p_jq_j^4+p_j^4q_j)+O(u^5)\right)\\
&=1+\frac{u^2}{2}\sum_{j=1}^{n/2}p_jq_j +
\frac{u^3}{3!}\sum_{j=1}^{n/2}p_jq_j(q_j-p_j)\\ & \qquad
+u^4\left[\frac{1}{4}\sum_{j_1\neq j_2}
(p_{j_1}q_{j_1})(p_{j_2}q_{j_2})+\frac{1}{4!}\sum_{j=1}^{n/2}(p_jq_j^4+p_j^4q_j
)\right]+O(u^5).
\end{split}\end{equation*}

\noindent In particular,

\begin{equation*}
\text{Var}[M_i]=E\left[\overline{M}_i^2\right]=\sum_{j=1}^{n/2}p_jq_j,
\end{equation*}

\noindent and consequently

\begin{equation*} \begin{split}
E\left[\overline{M}_i^4\right] &=6\sum_{j_1\neq j_2}(p_{j_1}q_{j_1})(p_{j_2}q_{j_2})+
\sum_{j=1}^{n/2}(p_jq_j^4+p_j^4q_j)\\
&\le 6\text{Var}^2[M_i]+\text{Var}[M_i].
\end{split} \end{equation*}

\noindent It is left to show that $\text{Var}[M_i]=O(n)$. Extending
the computation in (\ref{E[M]}),

\begin{equation*}\begin{split}
E[M_i(M_i-1)]&=\frac{n}{2}\(\frac{n}{2}-1\)\frac{\binom{n-2}{n/2-2}}
{\binom{n}{n/2}}\\ &=\frac{n^2}{16}+O(n) \\ &=E^2[M_i]+O(n).
\end{split}\end{equation*}

\noindent Therefore

\begin{equation*}
\text{Var}[M_i]=E[M_i(M_i-1)]+E[M_i]-E^2[M_i]=O(n).
\end{equation*}

\noindent This completes the proof of Theorem \ref{ThmSBr} (upper
bound). \epf

\bigskip

\begin{center}
{\Large {\bf Proof of Theorem \ref{ThmSBr}, lower bound.}}
\end{center}

We will show that for each $\epsilon>0$,

\[
P\( \pi \le \sigma\)=\Omega\(\( \alpha - \epsilon \)^n\),
\]

\noindent where

\[
\alpha = \sqrt[11]{ \frac{25497938851324213335}{22!} }=0.70879\dots.
\]

\noindent Introduce $\pi^*$ ($\sigma^*$ resp.), the permutation $\pi$
($\sigma$ resp.) with the element $n$ deleted. More generally, for $k
\le n$, $\pi^{k*}$ ($\sigma^{k*}$ resp.)  is the permutation of
$[n-k]$ obtained from $\pi$ ($\sigma$ resp.) by deletion of the $k$
largest elements, $n$, $n-1,\dots,n-k+1$. The key to the proof is the
following

\bigskip

\begin{Lem} \label{topkrowsLem} Let $k \in [n]$. If every northeastern 
submatrix of $M(\pi,\sigma)$ with at most $k$ rows contains at least as many
crosses as balls, and $\pi^{k*} \le \sigma^{k*}$, then $\pi \le \sigma$.
\end{Lem}

\bigskip

\noindent Before proceeding with the proof, we introduce one more bit
of notation. Let $G_n$ be the empty $n \times n$ grid depicted in the
$M(\cdot)$'s of Figure 2 above. Figure 5 below is a depiction of $G_5$
and an emboldened northeastern-corner $3\times 4$ subgrid of it,
denoted by $C$.

\bigskip
\begin{center}
\includegraphics[scale=0.7]{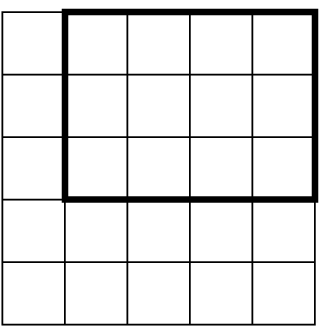} \\ {\small {\bf Figure 5:} $G_5$
and an emboldened subgrid $C$.}
\end{center}
\bigskip

\noindent If $C$ is any subgrid of $G_n$, then $M(\cdot \, | \, C)$
denotes the submatrix of $M(\cdot)$ that ``sits'' on $C$. To repeat,
the $(0,1)$-matrix criterion says that $\pi \le \sigma$ if and only if
for each northeastern-corner subgrid $C$ of $G_n$, we have at least as
many crosses as balls in $M(\pi,\sigma \, | \, C)$.

\bigskip

\bpf By the assumption, it suffices to show that the balls do not outnumber
crosses in $M(\pi,\sigma\,|\,C)$ for every such subgrid $C$ with
strictly more than $k$ rows.  Consider any such $C$.  Let $C^{(k)}$
denote the subgrid formed by the top $k$ rows of $C$. Given a
submatrix $A$ of $M(\pi)$ (of $M(\sigma)$ resp.), let $|A|$ denote the
number of columns in $A$ with a cross (a ball resp.). We need to show

\[
|M(\pi \, | \, C)| \ge |M(\sigma \, | \, C)|.
\]

\noindent By the assumption, we have $|M(\pi \, | \, C^{(k)})| \ge
|M(\sigma \, | \, C^{(k)})|$. Write $|M(\pi \, | \, C^{(k)})| =
|M(\sigma \, | \, C^{(k)})|+\lambda$, $\lambda\ge 0$. We now delete
the top $k$ rows from $M(\pi),M(\sigma)$ together with the $k$ columns
that contain the top $k$ crosses in the case of $M(\pi)$ and the $k$
columns that contain the top $k$ balls in the case of
$M(\sigma)$. This produces the matrices $M(\pi^{k*})$ and
$M(\sigma^{k*})$. In either case, we obtain the grid $G_{n-k}$
together with a new northeastern subgrid: $C(\pi^{k*})$ in the case of
$M(\pi)$ and $C(\sigma^{k*})$ in the case of $M(\sigma)$. Figure 6
below is a graphical visualization of this deletion process in the
special case $\pi=12534$, $\sigma=45132$, $k=2$ and $C$ the $3 \times
4$ northeastern subgrid of $G_5$. We have emboldened $C$ in
$M(\pi),M(\sigma)$, and the resulting $C(\pi^{2*}),C(\sigma^{2*})$ in
$M(\pi^{2*}),M(\sigma^{2*})$ respectively.

\bigskip
\begin{center}
\includegraphics[scale=0.7]{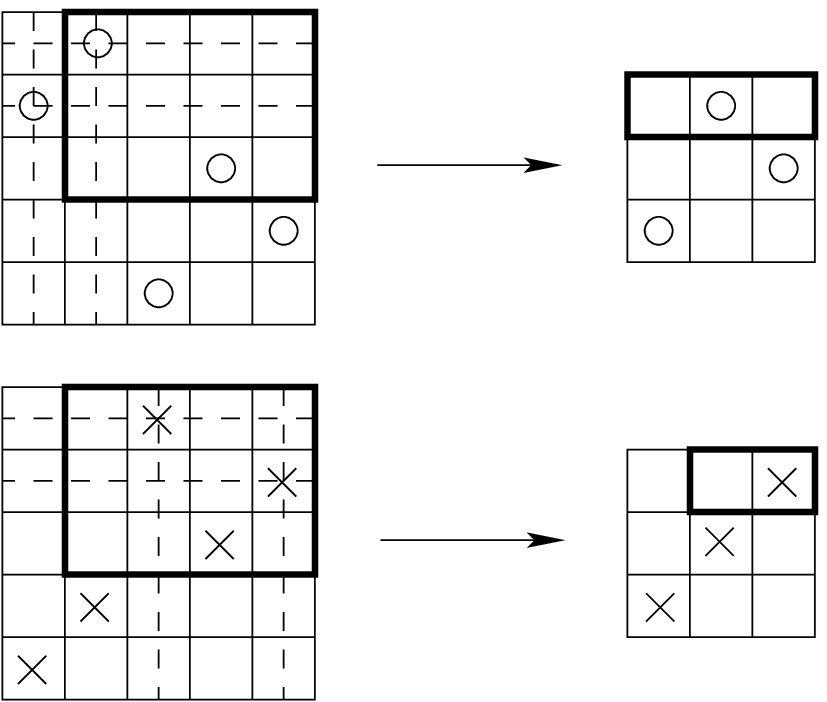} \\ {\small {\bf Figure
6:} Deletion of 2 largest elements of $\pi,\sigma$, and its affect on
$C$.}
\end{center}
\bigskip

\noindent Since we delete more columns in the case of $\pi$ than
$\sigma$, note that $C(\pi^{k*}) \subseteq C(\sigma^{k*})$ as
northeastern subgrids of $G_{n-k}$. In fact, these grids have the same
number of rows, but $C(\pi^{k*})$ has $\lambda$ fewer columns. Hence,
as $\pi^{k*} \le \sigma^{k*}$, we have

\[
|M(\pi^{k*} \, | \, C(\pi^{k*}))| \ge |M(\sigma^{k*} \, | \, C(\pi^{k*}))| \ge
 |M(\sigma^{k*} \, | \, C(\sigma^{k*}))|-\lambda.
\]

\noindent So

\begin{equation*} \begin{split}
|M(\pi \, | \, C)|&= |M(\pi \, | \, C^{(k)})|+|M(\pi^{k*} \, | \,
C(\pi^{k*}))| \\ &= |M(\sigma \, | \, C^{(k)})|+\lambda+|M(\pi^{k*} \, | \,
C(\pi^{k*}))| \\ &\ge |M(\sigma \, | \, C^{(k)})|+|M(\sigma^{k*} \, | \,
C(\sigma^{k*}))| \\ &= |M(\sigma \, | \, C)|,
\end{split} \end{equation*}

\noindent which proves the lemma. \epf

\bigskip

\noindent For each $k \le n$, let ${\cal E}_{n,k}$ denote the event
``every northeast submatrix of the top $k$ rows has at least as many crosses
as balls''. Then by Lemma \ref{topkrowsLem},

\[
\left\{ \pi \le \sigma \right\} \supseteq {\cal E}_{n,k} \cap \left\{
\pi^{k*} \le \sigma^{k*} \right\}.
\]

\noindent Now the events ${\cal E}_{n,k}$ and $\left\{\pi^{k*}\le
\sigma^{k*}\right\}$ are independent! So we get

\begin{equation} \label{LBest}
P\(\pi \le \sigma \) \ge P\({\cal E}_{n,k}\) P\(\pi^{k*} \le
\sigma^{k*} \).
\end{equation}

\noindent For the permutation $\pi$ ($\sigma$ resp.) introduce
$\ell_i(\pi)=\pi^{-1}(i)$ ($\ell_i(\sigma)=\sigma^{-1}(i)$ resp.), the
index of a column that contains a cross (a ball resp.)  at the
intersection with row $i$. In terms of the $\ell_i(\cdot)$'s, ${\cal
E}_{n,k}$ is the event: for each integer $j \le k$ and $m \le n$, the
number of $\ell_n(\pi),\ell_{n-1}(\pi),\dots,\ell_{n-j+1}(\pi)$ that
are $m$ at least is more than or equal to the number of
$\ell_n(\sigma),\ell_{n-1}(\sigma),\dots,\ell_{n-j+1}(\sigma)$ that are
$m$ at least. We could have replaced an integer $m\le n$ with a real number, 
which means that

$$
{\cal E}_{n,k}=\{(\pi,\sigma):\,(\boldsymbol{\ell}(\pi),\boldsymbol{\ell}
(\sigma))\in {\cal C}_k\},
$$ 

\noindent for some cone-shaped (Borel) set ${\cal C}_k \subset {\mathbb R}^{2k}$;
here $\boldsymbol{\ell}(\pi)=\{\ell_{n-i+1}(\pi)\}_{1\le i\le k}$,
$\boldsymbol{\ell}(\sigma)=\{\ell_{n-i+1}(\sigma)\}_{1\le i\le k}$.

Our task is to estimate sharply $P\({\cal E}_{n,k}\)$ for a fixed $k$,
and $n\to\infty$.  Observe first that $\boldsymbol{\ell}(\pi)$ and 
$\boldsymbol{\ell}(\sigma)$ are independent, and each uniformly distributed. 
For instance

$$ 
P\(\ell_n(\pi)=j_1,\dots, \ell_{n-k+1}(\pi)=j_k\)=\frac{1}{(n)_k},\quad
1\le j_1\neq\cdots\neq j_k\le n,
$$

\noindent where $(n)_k=n(n-1)\cdots(n-k+1)$. Since $(n)_k\sim n^k$ as
$n\to\infty$, $\ell_n(\pi), \dots,\ell_{n-k+1}(\pi)$ are almost independent
$[n]$-uniforms for large $n$, and fixed $k$.  Let us make this
asymptotic reduction rigorous.  Let $U$ be a uniform-$[0,1]$ random
variable, and let $U_1,\dots,U_n$ be independent copies of $U$. Then
each $\lceil nU_i \rceil$ is uniform on $[n]$, and it is easy to show
that

$$ 
P\(\lceil nU_1\rceil=i_1,\dots, \lceil nU_k\rceil =i_k \; | \;
\lceil nU_1\rceil\neq\cdots\neq \lceil nU_k\rceil\)=\frac{1}{(n)_k}.
$$

\noindent In other words, $\{\ell_{n-i+1}(\pi)\}_{1\le i\le k}$ has the same
distribution as $\lceil n\mathbf U\rceil:=\{\lceil nU_i\rceil \}_{1\le
i\le k}$ conditioned on the event ${\cal A}_{n,k}=\{\lceil nU_1\rceil\neq
\dots\neq\lceil nU_k\rceil\}$.  Analogously $\{\ell_{n-i+1}(\sigma)\}_{1\le
i\le k}$ is distributed as $\lceil n\mathbf V\rceil:=\{\lceil
nV_i\rceil\}_{1\le i\le k}$ conditioned on ${\cal B}_{n,k}=\{\lceil
nV_1\rceil\neq\cdots\neq \lceil nV_k\rceil\}$, where $V_1,\dots,V_k$
are independent $[0,1]$-uniforms, independent of $U_1,\dots, U_k$. We
will need yet another event ${\cal D}_{n,k}$ on which

$$ 
\min\{\min_{i\neq j}|U_i-U_j|,\,\min_{i\neq
j}|V_i-V_j|,\,\min_{i,j}|U_i-V_j|\}>1/n.
$$ 

\noindent Clearly on ${\cal D}_{n,k}$

$$
(\lceil n\mathbf U\rceil, \lceil n\mathbf V\rceil)\in
{\cal C}_k \Longleftrightarrow (\mathbf U,\mathbf V)\in {\cal C}_k;
$$ 

\noindent here $\mathbf U:=\{U_i\}_{1 \le i\le k}$, $\mathbf V:=\{V_i\}_{1 \le i\le k}$. 
In addition ${\cal D}_{n,k}\subseteq {\cal A}_{n,k}\cap {\cal B}_{n,k}$, and

$$ 
P({\cal D}_{n,k}^c)\le 2k^2P(|U_1-U_2|\le 1/n)\le 4k^2/n.
$$

\noindent Therefore

\begin{equation*} \begin{split}
P({\cal E}_{n,k})&=P((\boldsymbol{\ell}(\pi),\boldsymbol{\ell}(\sigma))\in {\cal C}_k)\\
&=\frac{P\left(\{(\lceil n\mathbf U\rceil,\lceil n\mathbf V\rceil)\in {\cal C}_k\}\cap\{{\cal A}_{n,k}
\cap {\cal B}_{n,k}\}\right)}{P({\cal A}_{n,k}\cap {\cal B}_{n,k})}\\
&=\frac{P\left(\{(\lceil n\mathbf U\rceil,\lceil n\mathbf V \rceil)\in {\cal C}_k\}\cap {\cal D}_{n,k}\right)+
O(P({\cal D}_{n,k}^c))}{1-O(P({\cal D}_{n,k}^c))}\\
&=\frac{P\left((\mathbf U,\mathbf V)\in {\cal C}_k\right)+O(k^2/n)}{1-O(k^2/n)}\\
&=Q_k+O(k^2/n),
\end{split} \end{equation*}

\noindent where $Q_k=P\left((\mathbf U,\mathbf V)\in {\cal C}_k\right)$. Let us write $P_n=P(\pi \le \sigma)$. 
Using (\ref{LBest}) and the last estimate, we obtain then

\[
P_n\ge Q_kP_{n-k}\(1+O(k^2/n)\)=Q_kP_{n-k}\exp\(O(k^2/n)\),\quad n>k.
\]

\noindent Iterating this inequality $\lfloor n/k\rfloor$ times gives

\[
P_n\ge Q_k^{\lfloor n/k\rfloor}P_{n-\lfloor n/k\rfloor k}\exp\left[\sum_{j=0}^{\lfloor n/k\rfloor -1}
O\left(\frac{k^2}{n-jk}\right)\right].
\]

\noindent Since the sum in the exponent is of order $O(k^2\log n)$, we get

\[
\liminf \sqrt[n]{P_n} \ge \sqrt[k]{Q_k},\quad\forall k\ge 1.
\]

\noindent Thus

\[
\liminf \sqrt[n]{P_n} \ge \sup_k \sqrt[k]{Q_k}.
\]

\noindent Therefore, for each $k$ and $\epsilon\in \(0,\sqrt[k]{Q_k}\)$, we have

\begin{equation} \label{LBorder}
P_n=\Omega\( \( \sqrt[k]{Q_k}-\epsilon \)^n \).
\end{equation}

\noindent Next

\bigskip

\begin{Lem} \label{Q(k)supMult} As a function of $k$, $Q_k$ is supermultiplicative,
i.e. $Q_{k_1+k_2}\ge Q_{k_1}Q_{k_2}$ for all $k_1,k_2\ge 1$. Consequently
there exists $\lim_{k \to \infty} \sqrt[k]{Q_k}$, and moreover

\[
\lim_{k \to \infty} \sqrt[k]{Q_k}=\sup_{k\geq 1} \sqrt[k]{Q_k}.
\] 

\end{Lem}

\bigskip

\noindent Thus we expect that our lower bound would probably improve as $k$ increases.  

\bigskip

\bpf $Q_k$ is the probability of the event $E_k=\{(\mathbf U^{(k)},\mathbf V^{(k)})\in {\cal C}_k\}$;
here $\mathbf U^{(k)}:=\{U_i\}_{1\le i\le k}$, $\mathbf V^{(k)}:=\{V_i\}_{1\le i\le k}$.
Explicitly, for each $j\le k$ and each $c\in [0,1]$, the number of $U_1,
\dots, U_j$ not exceeding $c$ is at most the number of $V_1,\dots,V_j$ not
exceeding $c$. So $Q_{k_1+k_2}=P(E_{k_1+k_2})$, $Q_{k_1}=P(E_{k_1})$, while
$Q_{k_2}=P(E_{k_2})=P(E_{k_2}^*)$. Here the event $E_{k_2}^*$ means that for each $j^* \le
k_2$ and each $c \in [0,1]$, the number of $U_i$, $i=k_1+1,\dots,k_1+j^*$, not exceeding $c$ 
is at most the number of $V_i$, $i=k_1+1,\dots,k_1+j^*$, not exceeding $c$. The events
$E_{k_1}$ and $E_{k_2}^*$ are independent. Consider the
intersection of $E_{k_1}$ and $E_{k_2}^*$. There are two
cases:

\begin{enumerate}
\item[1)] $j \le k_1$. Then the number of $U_i$, $i \le j$ not exceeding
$c$ is at most the number of $V_i$, $i \le j$ not exceeding
$c$, as $E_{k_1}$ holds.

\item[2)] $k_1<j \le k_1+k_2$. Then the number of $U_i$, $i \le j$,
not exceeding $c$ is at most the number of $V_i$, $i \le k_1$ not
exceeding $c$ (as $E_{k_1}$ holds), plus the number of $V_i$,
$k_1< i \le j$, not exceeding $c$ (as $E_{k_2}^*$ holds).
The total number of these $V_i$ is the number of all
$V_i$, $i \le j$, that are at most $c$, $c \in [0,1]$.
\end{enumerate}

\noindent So $E_{k_1+k_2} \supseteq E_{k_1}\cap E_{k_2}^*$,
and we get $Q_{k_1+k_2}\ge Q_{k_1}Q_{k_2}$. The rest of the statement follows 
from a well-known result about super(sub)multiplicative sequences 
(see P\'olya and Szeg\"o \cite{PS}). \epf

\bigskip

Given $1 \le j \le i \le k$, let $U_{i,j}$ ($V_{i,j}$ resp.) denote the $j$-th element in the increasing rearrangement of $U_1,\dots,U_i$ ($V_1,\dots,V_i$ resp.). Then, to put it another way, $Q_k$ is the probability that the $k$ Ehresmann conditions are met by the independent $k$-dimensional random vectors $\bf U$ and $\bf V$, both of which have independent entries. That is, we check $U_{i,j}>V_{i,j}$ for each $1 \le j \le i \le k$ by performing element-wise comparisons in the following tableaux:

\[
\begin{array}{c}
\begin{array}{ccccc}
U_{k,1} & U_{k,2} & U_{k,3} & \cdots & U_{k,k} \\ \vdots & \vdots & \vdots & \reflectbox{$\ddots$} & \\ U_{3,1} & U_{3,2} & U_{3,3} & & \\ U_{2,1} & U_{2,2} & & & \\ U_{1,1} & & & &
\end{array}
\qquad \qquad
\begin{array}{ccccc}
V_{k,1} & V_{k,2} & V_{k,3} & \cdots & V_{k,k} \\ \vdots & \vdots & \vdots & \reflectbox{$\ddots$} & \\ V_{3,1} & V_{3,2} & V_{3,3} & & \\ V_{2,1} & V_{2,2} & & & \\ V_{1,1} & & & &
\end{array}
\end{array}.
\]

\noindent What's left is to explain how we determined $\alpha=0.70879...$.

It should be clear that whether or not $(\mathbf U^{(k)},\mathbf V^{(k)})$
is in ${\cal C}_k$ depends only on the size ordering of $U_1,\dots, U_k,V_1,\dots,
V_k$. There are $(2k)!$ possible orderings, all being equally likely.
Thus $Q_k=N_k/(2k)!$. Since the best constant in the lower exponential
bound is probably $\lim_{k\to\infty}\sqrt[k]{Q_k}$, our task was to compute $N_k$ 
for $k$ as large as our computer could handle. (``Probably'', because we do not know for certain 
that $\sqrt[k]{Q_k}$ increases with $k$.)

Here is how $N_k$ was tabulated. Recursively, suppose we have determined all $N_{k-1}$
orderings of $x_1$, $\dots$, $x_{k-1}$, $y_1$, $\dots$, $y_{k-1}$ such that $(\mathbf x^{(k-1)},
\mathbf y^{(k-1)}) \in {\cal C}_{k-1}$. Each such ordering can be assigned a $2(k-1)$-long sequence of 
$0$'s and $1$'s, $0$'s for $x_i$'s and $1$'s for $y_j$'s, $1\le i,j\le k-1$. Each such sequence 
meets the ballot-theorem condition: as we read it from {\it left to right\/} the number of $1$'s 
never falls below the number of $0$'s. We also record the {\it multiplicity\/} of each sequence, which
is the number of times it is encountered in the list of all $N_{k-1}$ orderings. The knowledge of all 
$2(k-1)$-long ballot-sequences together with their multiplicities is all we need to compile
the list of all $2k$-long ballot-sequences with their respective multiplicities. 

For $k=1$, there is only
one ballot-sequence to consider, namely $10$, and its multiplicity is $1$. So $N_1=1$,
and 

\[
Q_1=1/2!.
\]

\noindent Passing to $k=2$, we
must count the number of ways to insert $1$ and $0$ into $10$ so that we
get a 4-long ballot-sequence of two $0$'s and two $1$'s. Inserting $1$ at
the beginning, giving $110$, we can insert $0$ into positions $2$, $3$ or $4$,
producing three ballot-sequences

\[
1010, \quad 1100, \quad 1100,
\]

\noindent respectively. (Inserting $0$ into position $1$ would have resulted in $0110$ which
is not a ballot-sequence.) Similarly, inserting $1$ into position $2$, we get
$110$, and inserting $0$ under the ballot condition gives three ballot-sequences

\[
1010, \quad 1100, \quad 1100.
\]

\noindent Finally, inserting $1$ at the end, giving $101$, we can only
insert $0$ at the end, obtaining one ballot-sequence

\[
1010.
\]

\noindent Hence, starting from the ballot-sequence $10$ of multiplicity $1$,
we have obtained two $4$-long ballot-sequences, $1010$ of multiplicity $3$ and $1100$ of
multiplicity $4$. Therefore $N_2=3+4=7$, and

\[
Q_2=7/4!.
\]

\noindent Pass to $k=3$. Sequentially we insert $1$ in each of $5$ positions in
the ballot-sequence $1010$, and then determine all positions for the new $0$ which
would result in a $6$-long ballot-sequence. While doing this we keep track of how many
times each $6$-long ballot-sequence is encountered. Multiplying these numbers by $3$, the 
multiplicity of $1010$, we obtain a list of $6$-long ballot-sequences spawned by $1010$ with the 
number of their occurrences. We do the same with the second sequence $1100$. Adding the numbers of 
occurrences of each $6$-long ballot-sequence for $1010$ and $1100$, we arrive at the following list
of five $6$-long ballot-sequences with their respective multiplicities:

\begin{equation*} \begin{split}
111000 &: 36, \\ 110100 &: 32, \\ 110010 &: 24, \\ 101100 &: 24, \\
101010 &: 19.
\end{split} \end{equation*} 

\noindent Therefore $N_3=36+32+24+24+19=135$, and 

\[
Q_3=135/6!.
\]

\noindent The younger coauthor wrote a computer program for this algorithm. Pushed
to its limit, the computer delivered the following table:

\bigskip
\begin{center}
\begin{tabular}{r|r|l|l}
\multicolumn{1}{c|}{$k$} & \multicolumn{1}{c|}{$N_k=(2k)!Q_k$} &
\multicolumn{1}{c|}{$Q_k=N_k/(2k)!$} &
\multicolumn{1}{c}{$\sqrt[k]{Q_k}$} \\ \hline 1 & 1 & $0.50000\dots$ &
$0.50000\dots$ \\ 2 & 7 & $0.29166\dots$ & $0.54006\dots$ \\ 3 & 135 &
$0.18750\dots$ & $0.57235\dots$ \\ 4 & 5193 & $0.12879\dots$ &
$0.59906\dots$ \\ 5 & 336825 & $0.09281\dots$ & $0.62162\dots$ \\ 6 &
33229775 & $0.06937\dots$ & $0.64101\dots$ \\ 7 & 4651153871 &
$0.05335\dots$ & $0.65790\dots$ \\ 8 & 878527273745 & $0.04198\dots$ &
$0.67280\dots$ \\ 9 & 215641280371953 & $0.03368\dots$ &
$0.68608\dots$ \\ 10 & 66791817776602071 & $0.02745\dots$ &
$0.69800\dots$ \\ 11 & 25497938851324213335 & $0.02268\dots$ &
$0.70879\dots$ \\
\end{tabular}
\end{center}
\bigskip

\noindent Combining (\ref{LBorder}) and the value of
$\sqrt[11]{Q_{11}}$ in this table, we see that for each $\epsilon>0$,

\[
P_n= \Omega\( \( \sqrt[11]{Q_{11}} - \epsilon \)^n \)=\Omega\((0.708...-\epsilon)^n\). \epf
\]

\bigskip

\noindent The numbers $\sqrt[k]{Q_k}$ increase steadily for $k<12$,
so at this moment we would not rule out the tantalizing possibility that
$\sqrt[k]{Q_k}\to 1$ as $k\to\infty$. Determination of the actual
limit is a challenging open problem. The proof just given only involves mention of 
the $(0,1)$-matrix criterion, but it was the Ehresmann criterion that actually inspired 
our initial insights.

\bigskip

\begin{center}
{\Large {\bf Weak Bruhat Order, preliminaries.}}
\end{center}

Recall that  $\pi$ precedes $\sigma$ in the
\textit{weak} Bruhat order ($\pi \preceq\sigma$) if and only if 
there is a chain $\sigma=\omega_1\to\cdots\to\omega_s=\pi$ where
each $\omega_t$ is a simple reduction of $\omega_{t-1}$, i.e. obtained 
by transposing two adjacent elements $\omega_{t-1}(i)$, $\omega_{t-1}(i+1)$
such that $\omega_{t-1}(i)>\omega_{t-1}(i+1)$. Clearly the weak Bruhat order is more
restrictive, so that $\pi\preceq \sigma$ impies $\pi\le\sigma$.  In
particular, $P(\pi\preceq\sigma)\le P(\pi\le\sigma)$, hence (Theorem \ref{ThmSBr})
$P(\pi\preceq \sigma)=O(n^{-2})$. We will show that, in fact, this probability
is exponentially small. The proof is based on an inversion set criterion for $\pi\preceq
\sigma$ implicit in \cite[pp. 135-139]{Be}.

\bigskip

\begin{Lem} \label{LemWBrCrit} Given $\omega\in S_n$, introduce
the set of non-inversions of $\omega$:

\[
E \( \omega \):= \left\{ \( i,j \) \; : \; i<j,\, \omega^{-1} \( i\) <
\omega^{-1} \( j\) \right\}.
\]

\noindent $\pi \preceq \sigma$ if and only if $E \( \pi
\)\supseteq E \( \sigma \)$.
\end{Lem}

\bpf Assume $\pi \preceq \sigma$. Then there exists a chain of simple 
reductions $\omega_t$, $1\le t\le s$, connecting $\sigma=\omega_1$ and 
$\pi=\omega_s$. By the definition of a simple reduction, for each $t>1$
there is $i=i_t<n$ such that $E(\omega_t)=E(\omega_{t-1})\cup 
\{\(\omega_t(i),\omega_t(i+1)\)\}$, where $\omega_t(i)=\omega_{t-1}(i+1)$,
$\omega_t(i+1)=\omega_{t-1}(i)$, and $\omega_{t-1}(i)>\omega_{t-1}(i+1)$.
So the set $E(\omega_t)$ increases with $t$, hence $E(\pi)\supseteq
E(\sigma)$.

Conversely, suppose $E\(\pi\) \supseteq E\(\sigma\)$. Since a
permutation $\omega$ is uniquely determined by its $E(\omega)$, we may
assume $E\(\pi\) \supsetneq E\(\sigma\)$. 

\bigskip

\noindent {\bf Claim~} If $E\(\pi\) \supsetneq E\(\sigma\)$, then
there exists $u<v\le n$ such that $\(v,u\)$ is an adjacent inversion
of $\sigma$, but $\(u,v\)\in E\(\pi\)$.

\bigskip

\noindent Assuming validity of the claim, we ascertain existence
of an adjacent inversion $\(v,u\)$ in $\sigma$ with $\(u,v\)\in E\(\pi\)$. 
Interchanging the adjacent elements $u$ and $v$ in $\sigma=\omega_1$, we
obtain a simple reduction $\omega_2$, with $E(\omega_1)\subset E(\omega_2)
\subseteq E(\pi)$. If $E(\omega_2)=E(\pi)$ then $\omega_2=\pi$, and we stop.
Otherwise we determine $\omega_3$, a simple reduction of $\omega_2$,
with $E(\omega_2)\subset E(\omega_3)\subseteq E(\pi)$ and so on. Eventually
we determine a chain of simple reductions connecting $\sigma$ and $\pi$,
which proves that $\pi\preceq \sigma$. 

\bigskip

\noindent {\sc Proof of Claim.~} The
claim is obvious for $n=1,2$. Assume inductively that the
claim holds for permutations of length $n-1 \ge 2$. Let $\pi,\sigma \in S_n$
and $E(\pi)\supsetneq E(\sigma)$. As in the proof of Theorem \ref{ThmSBr},
let $\ell_n(\pi)=\pi^{-1}(n)$, $\ell_n(\sigma)=\sigma^{-1}(n)$, and $\pi^*$,
$\sigma^*$ are obtained by deletion of $n$ from $\pi$ and $\sigma$.
Since $E\(\pi\) \supseteq E\(\sigma\)$, we have $E\(\pi^*\)
\supseteq E\(\sigma^*\)$. Suppose first that $E\(\pi^*\) =
E\(\sigma^*\)$. Then $\pi^*=\sigma^*$, and as $E\(\pi\) \supsetneq
E\(\sigma\)$, we must have $\ell_n(\pi) > \ell_n(\sigma)$, i.e.
$\ell_n(\sigma)<n$. Setting
$v=n$ and $u=\sigma\(\ell_n(\sigma)+1\)$, we obtain an
adjacent inversion $\(v,u\)$ in $\sigma$ with $\(u,v\)\in E\(\pi\)$.

Alternatively, $E\(\pi^*\) \supsetneq E\(\sigma^*\)$. By inductive
hypothesis, there exists $u<v\le n-1$ such that $\(v,u\)$ is an
adjacent inversion of $\sigma^*$, but $\(u,v\) \in E\(\pi^*\)$. Now
insert $n$ back into $\pi^*,\sigma^*$, recovering $\pi$ and
$\sigma$. If $n$ sits to the right of $u$ or to the left of $v$ in
$\sigma$, then $\(v,u\)$ is still an adjacent inversion of $\sigma$.
Otherwise $n$ is sandwiched between $v$ on the left and $u$ on the
right. Therefore $\(n,u\)$ is an adjacent inversion in $\sigma$.
On the other hand $\(v,n\)\in E\(\sigma\)$, so since
$E\(\pi\) \supseteq E\(\sigma\)$, we have $\(v,n\) \in E\(\pi\)$
also. Hence, the triple $\(u,v,n\)$ are in exactly this order (not
necessarily adjacent) in $\pi$. Therefore the adjacent inversion
$\(n,u\)$ in $\sigma$ is such that $\(u,n\) \in E\(\pi\)$, and 
this proves the inductive step. \epf

\bigskip

\noindent  Denote by $\bar{\omega}$ the
permutation $\omega$ reversed in rank. For example, with $\omega=13254$ 
we have $\bar{\omega}=53412$. Then it is easy to see that

\[
E \( \pi \)\supseteq E \( \sigma \) \Longleftrightarrow E \( \bar{\pi}
\)\subseteq E \( \bar{\sigma} \).
\]

\noindent By Lemma \ref{LemWBrCrit}, these statements are equivalent
to

\[
\pi \preceq \sigma \Longleftrightarrow  \bar{\sigma} \preceq \bar{\pi}.
\]

\noindent We immediately obtain the following corollary to Lemma
\ref{LemWBrCrit}:

\bigskip

\begin{Cor} \label{CorEi} For $\omega\in S_n$, define

\[
E_i\(\omega\):= \left\{ j<i \, : \, \(j,i\) \in E\(\omega\) \right\}, \quad
1 \leq i \leq n.
\]

\noindent Then

\[
E \(\omega\)= \bigsqcup_{i=1}^n \left\{ \(j,i\) \, : \, j \in E_i\(\omega\)
\right\},
\]

\noindent and consequently

\[
\pi \preceq \sigma \Longleftrightarrow E\(\pi\) \supseteq E\(\sigma\)
\Longleftrightarrow E_i\(\pi\) \supseteq E_i\(\sigma\), \;\; \forall i
\leq n \textrm{.} \epf
\]

\end{Cor}

\bigskip

\noindent Next,

\bigskip

\begin{Lem} \label {LemSubMult} For $\pi,\sigma \in S_n$ selected independently and uniformly 
at random, let $P_n^*:=P\( \pi \preceq \sigma \)$. $P_n^*$ is submultiplicative,
i.e. for all $n_1,n_2\ge 1$

\[
P_{n_1+n_2}^* \leq P_{n_1}^* P_{n_2}^* \textrm{.}
\]

\noindent Consequently there exists $\lim_{n\to\infty}\sqrt[n]{P_n^*}=\inf_{n\ge 1}\sqrt[n]{P_n^*}$.
\end{Lem}

\bpf Let $\pi,\sigma$ be two permutations of
$\left[n_1+n_2\right]$. Then $\pi \preceq \sigma$ if and only if

\[
E_i \(\pi\)\supseteq E_i \(\sigma\), \quad 1\leq i \leq n_1+n_2
\textrm{.}
\]

\noindent Using these conditions for $i\leq n_1$, we see that

\[
\pi\left[1,2,\dots,n_1\right] \preceq \sigma\left[1,2,\dots,n_1\right]
\textrm{.}
\]

\noindent Here $\pi\left[1,2,\dots,n_1\right]$, say, is what is left of
the permutation $\pi$ when the elements $n_1+1,\dots,n_1+n_2$ are
deleted.

Likewise, $\pi \preceq \sigma$ if and only if

\[
E_i\( \bar{\pi}\)\subseteq E_i\( \bar{\sigma}\), \quad 1\leq i\leq
n_1+n_2 \textrm{.}
\]

\noindent Using these conditions for $i \leq n_2$, we see that

\[
\pi\left[n_1+1,\dots,n_1+n_2\right] \preceq
\sigma\left[n_1+1,\dots,n_1+n_2\right] \textrm{.}
\]

\noindent Now, since $\pi$ and $\sigma$ are uniformly random and
mutually independent, so are the four permutations

\[
\pi\left[1,\dots,n_1\right], \quad \pi\left[n_1+1,\dots,n_1+n_2\right],
\quad \sigma\left[1,\dots,n_1\right], \quad
\sigma\left[n_1+1,\dots,n_1+n_2\right] \textrm{.}
\]

\noindent Hence,

\[
P\( \pi \preceq \sigma\)\leq P\(\pi\left[1,\dots,n_1\right] \preceq
                  \sigma\left[1,\dots,n_1\right]\) \cdot
                  P\(\pi\left[n_1+1,\dots,n_1+n_2\right] \preceq
                  \sigma\left[n_1+1,\dots,n_1+n_2\right]\) \textrm{,}
\]

\noindent so that

\[
P_{n_1+n_2}^* \leq P_{n_1}^* P_{n_2}^* \textrm{.} \epf
\]

\bigskip

\begin{center}
{\Large {\bf Proof of Theorem \ref{ThmWBr}, upper bound.}}
\end{center}

We will show that for each $\epsilon>0$,

\[
P_n^*=O((\beta + \epsilon)^n),
\]

\noindent where

\[
\beta=\sqrt[6]{\frac{1065317}{12!}}=0.36129\dots.
\]

\noindent The proof of this upper bound for $P_n^*$ parallels the proof
of the lower bound for $P_n$ in Theorem \ref{ThmSBr}. As in that proof, given 
$k\ge 1$, let $\pi^{k*}$ and $\sigma^{k*}$ be obtained 
by deletion of the elements $n,\dots,n-k+1$ from $\pi$ and $\sigma$, and
let $\ell_i(\pi)=\pi^{-1}(i)$, $\ell_i(\sigma)=\sigma^{-1}(i)$, $n-k+1\le i\le n$.  
In the notations of the proof of  Lemma \ref{LemSubMult}, $\pi^{k*}=\pi[1,\dots,n-k]$
and $\sigma^{k*}=\sigma[1,\dots,n-k]$, and we saw that $\pi^{k*}\preceq \sigma^{k*}$
if $\pi\preceq \sigma$. Our task is to find the conditions these $\ell_i(\cdot)$'s
must satisfy if $\pi\preceq \sigma$ holds.

To start, notice that

\[
\pi \preceq \sigma \Longrightarrow |E_n(\pi)| \ge |E_n(\sigma)|
\Longleftrightarrow \ell_n(\pi) \ge \ell_n(\sigma).
\]

\noindent Next

\[
\pi \preceq \sigma \Longrightarrow \pi^* \preceq \sigma^*
\Longrightarrow \ell_{n-1}(\pi) \ge \ell_{n-1}\(\pi^*\) \ge
\ell_{n-1}\(\sigma^*\) \ge \ell_{n-1}(\sigma)-1,
\]

\noindent as deletion of $n$ from $\pi,\sigma$ decreases
the location of $n-1$ in each permutation by at most one. In general,
for $0 < j < k$ we get

\[
\pi \preceq \sigma \Longrightarrow \pi^{j*} \preceq \sigma^{j*}
\Longrightarrow \ell_{n-j}(\pi) \ge \ell_{n-j}(\sigma)-j.
\]

\noindent  So, introducing $\boldsymbol{\ell}(\pi)=\{\ell_{n-i+1}(\pi)
\}_{1\le i\le k}$ and $\boldsymbol{\ell}(\sigma)=\{\ell_{n-i+1}(\sigma)\}_{1\le i\le k}$,

\begin{eqnarray} \label{piLEsigmaL(n,k)}
&\left\{ \pi \preceq \sigma \right\} \subseteq \{(\boldsymbol{\ell}(\pi),\boldsymbol{\ell}
(\sigma)\in {\cal S}_k\}, \\
&{\cal S}_k:= \left\{(\mathbf x,\mathbf y)\in \mathbb R^{2k}:\,  x_j\ge y_j-j+1,\,1\le j\le k\right\}. \nonumber
\end{eqnarray} 

\noindent  In addition, on $\{\pi\preceq\sigma\}$ every pair of elements,
which forms an inversion in $\pi$, also forms an inversion in $\sigma$.
Applying this to the elements $n-k+1,\dots,n$, we have then 

\begin{eqnarray} \label{piLEsigmaT(n,k)}
&\left\{ \pi \preceq \sigma \right\} \subseteq \{(\boldsymbol{\ell}(\pi),\boldsymbol{\ell}(\sigma))\in {\cal T}_k\},\\
&{\cal T}_k:=\{(\mathbf x,\mathbf y)\in \mathbb R^{2k}:\,\forall  1\le i<j\le k \, ,
x_i< x_j\Longrightarrow  y_i<y_j\}. \nonumber
\end{eqnarray}

\noindent Combining  (\ref{piLEsigmaL(n,k)}) and
(\ref{piLEsigmaT(n,k)}), we get

\[
\left\{ \pi \preceq \sigma \right\} \subseteq \{(\boldsymbol{\ell}(\pi),\boldsymbol{\ell}
(\sigma))\in {\cal S}_k\cap {\cal T}_k\}
\cap \left\{ \pi^{k*} \preceq \sigma^{k*} \right\}.
\]

\noindent So, since the two events on the right are independent,

\begin{equation} \label{Pn(*)UppEst}
P_n^* \le P( (\boldsymbol{\ell}(\pi),\boldsymbol{\ell}(\sigma))\in {\cal S}_k\cap {\cal T}_k)P_{n-k}^*.
\end{equation}

\noindent It remains to estimate $P( (\boldsymbol{\ell}(\pi),\boldsymbol{\ell}(\sigma))\in 
{\cal S}_k\cap {\cal T}_k)$.  As in the proof of Theorem \ref{ThmSBr} (lower bound), we observe that  $(\boldsymbol{\ell}(\pi),\boldsymbol{\ell}(\sigma))$
has the same distribution as $(\lceil n\mathbf U\rceil,\lceil n\mathbf V\rceil)$,
conditioned on 

$$
{\cal A}_{n,k}\cap {\cal B}_{n,k}=\{\lceil nU_1\rceil\neq\cdots\neq\lceil nU_k\rceil\}
\cap\{\lceil nV_1\rceil\neq\cdots\neq\lceil nV_k\rceil\}.
$$

\noindent Here $U_1,\dots,U_k,V_1,\dots,V_k$ are independent $[0,1]$-uniforms. Then

$$
P( (\boldsymbol{\ell}(\pi),\boldsymbol{\ell}(\sigma))\in {\cal S}_k\cap {\cal T}_k)=
\frac{P(\{(\lceil n\mathbf U\rceil,\lceil n\mathbf V\rceil)\in {\cal S}_k\cap 
{\cal T}_k\}\cap {\cal C}_{n,k})}{P({\cal C}_{n,k})},
\quad {\cal C}_{n,k}={\cal A}_{n,k}\cap {\cal B}_{n,k}.
$$

\noindent Introduce the event $\tilde{{\cal D}}_{n,k}$ on which

$$
\min\{\min_{i\neq j}|U_i-U_j|,\,\min_{i\neq j}|V_i-V_j|,\,\min_{i,j}|U_i-V_j|,\,k^{-1}\min_j|U_j-V_j|\}>
1/n.
$$

\noindent Certainly  $\tilde{{\cal D}}_{n,k}\subseteq {\cal C}_{n,k}$ and, thanks to
the factor $1/k$ by $\min_j|U_j-V_j|$, on $\tilde{{\cal D}}_{n,k}$

\begin{equation*}
\lceil nU_j\rceil\ge \lceil nV_j\rceil-j+1\Longrightarrow U_j\ge V_j-k/n\Longrightarrow
U_j> V_j.
\end{equation*}

\noindent Therefore, on $\tilde{{\cal D}}_{n,k}$,

\begin{eqnarray*}
&(\lceil n\mathbf U\rceil,\lceil n\mathbf V\rceil) \in {\cal S}_k\cap {\cal T}_k \Longrightarrow
(\mathbf U,\mathbf V)\in \tilde{{\cal S}}_k\cap {\cal T}_k,\\
&\tilde{{\cal S}}_k:=\{(\mathbf x,\mathbf y)\in \mathbb R^{2k}:\,x_j>y_j,\,1\le j\le k\}.
\end{eqnarray*}

\noindent Clearly $\tilde{{\cal S}}_k\cap {\cal T}_k$ is a cone-shaped subset of $\mathbb R^{2k}$.  
In addition, $P(\tilde{{\cal D}}_{n,k}^c)=O(k^2/n)$. Hence

\begin{equation*} \begin{split}
P( (\boldsymbol{\ell}(\pi),\boldsymbol{\ell}(\sigma))\in {\cal S}_k\cap {\cal T}_k)&
\le \frac{P((\mathbf U,\mathbf V)\in \tilde{{\cal S}}_k\cap {\cal T}_k)+O(P(\tilde{{\cal D}}_{n,k}^c))}
{1-O(P(\tilde{{\cal D}}_{n,k}^c))}\\
&=Q^*_k(1+O(k^2/n)),\quad Q_k^*:=P((\mathbf U,\mathbf V)\in \tilde{{\cal S}}_k\cap {\cal T}_k).
\end{split} \end{equation*}

\noindent This and (\ref{Pn(*)UppEst}) imply 

\[
P_n^* \le Q_k^*P_{n-k}^*\exp(O(k^2/n)).
\]

\noindent Hence, as in the proof of Theorem \ref{ThmSBr} (lower bound),

\[
\limsup \sqrt[n]{P_n^*} \le \sqrt[k]{Q_k^*},\quad\forall\, k\ge 1,
\]

\noindent and so

\begin{equation}\label{wUBPn*}
P_n^*=O((\sqrt[k]{Q_k^*}+\epsilon)^n),\quad k\ge 1,\,\epsilon>0.
\end{equation}

\noindent Furthermore, from the definition of $Q_k^*$, it follows directly that
$Q_k^*$ is {\it submultiplicative\/}, i.e. 

$$
Q_{k_1+k_2}^*\le Q_{k_1}^*Q_{k_2}^*,\quad k_1,k_2\ge 1.
$$

\noindent Therefore (\cite{PS} again)

$$
\lim_{k\to\infty}\sqrt[k]{Q_k^*}=\inf_{k\ge 1}\sqrt[k]{Q_k^*}.
$$

\noindent So the further we can push tabulation of $Q_k^*$,
the better our exponential upper bound for $P_n^*$ would
probably be. (``Probably'', because we do not have a proof that
$\sqrt[k]{Q_k^*}$ decreases with $k$.)

As in the case of $Q_k$, $Q_k^*=N_k^*/(2k)!$. Here, by the definition of
the sets $\tilde{{\cal S}}_k$ and ${\cal T}_k$, $N_k^*$ is the total
number of ways to order $x_1,\dots,x_k,y_1,\dots,y_k$ so that two
conditions are met: (1) for each $j$, $x_j$ is to the right
of $y_j$; (2) for all $i<j$, if $x_i$ is to the left of $x_j$ then
$y_i$ is to the left of $y_j$.

It is instructive first to evaluate $N_k^*$ by hand for $k=1,2$.
$N_1^*=1$ as there is only one sequence, $y_1x_1$, meeting the conditions
(1), (2). Passing to $N_2^*$, we must decide how to insert $y_2$ and
$x_2$ into the sequence $y_1x_1$ in compliance with conditions (1), (2).
First of all, $y_2$ has to precede $x_2$. If we insert $x_2$ at
the beginning of  $y_1x_1$, giving $x_2y_1x_1$, then
we can only insert $y_2$ at the beginning of this triple, giving

\[
y_2x_2y_1x_1.
\]

\noindent Alternatively, inserting $x_2$ in the middle of
$y_1x_1$, we have $2$ possibilities for insertion of $y_2$, and we
get two admissible orderings,

\[
y_2y_1x_2x_1, \qquad y_1y_2x_2x_1.
\]

\noindent Finally, insertion of $x_2$ at the end of $y_1x_1$ brings
the condition (2) into play as we now have $x_1$ preceding $x_2$, and so 
$y_1$ must precede $y_2$. Consequently, we get two admissible orderings, 

\[
y_1y_2x_1x_2, \qquad y_1x_1y_2x_2.
\]

\noindent Hence $N_2^*=1+2+2=5$. Easy so far! However, passing to
$k=3$ is considerably more time-consuming than it was for computation of $N_3$
in the proof of the lower bound in Theorem \ref{ThmSBr}. There, once we had determined
the $N_2$ admissible orderings, we could afford not to keep track of relative
orderings of $x_1,\dots,x_{k-1}$, and of $y_1,\dots,y_{k-1}$, whence the
coding by $1$'s and $0$'s. All we needed for passing from $k-1$ to $k$ was
the list of all binary ballot-sequences of length $2(k-1)$ together with
their multiplicities. Here the nature of the conditions (1), (2)
does not allow lumping various sequences together, and we have to preserve
the information of relative orderings of $x$'s, and relative orderings of
$y$'s. This substantial complication seriously inhibits the computer's
ability to compute $N_k^*$ for $k$ as large as in the case of $N_k$.

To get a feeling for how sharply the amount of computation increases for $k=3$,
let us consider one of the $N_2^*=5$
admissible sequences, namely $y_2x_2y_1x_1$. As above, we write down all possible
ways to insert $y_3$ and $x_3$ into this sequence so that (1) and (2)
hold. Doing this, we produce the $10$ sequences:

\begin{equation*} \begin{split}
y_3x_3y_2x_2y_1x_1, \qquad y_3y_2x_3x_2y_1x_1, \\
y_2y_3x_3x_2y_1x_1, \qquad y_2y_3x_2x_3y_1x_1, \\
y_2x_2y_3x_3y_1x_1, \qquad y_2y_3x_2y_1x_3x_1, \\
y_2x_2y_3y_1x_3x_1, \qquad y_2x_2y_1y_3x_3x_1, \\
y_2x_2y_1y_3x_1x_3, \qquad y_2x_2y_1x_1y_3x_3.
\end{split} \end{equation*} 

\noindent We treat similarly the other four sequences from the $k=2$ case, 
eventually arriving at $N_3^*=55$. We wouldn't even think of computing
$N_4^*$ by hand.

Once again the computer programming to the rescue! Here is the table produced 
by the computer after a substantial running time:

\bigskip
\begin{center}
\begin{tabular}{r|r|l|l}
\multicolumn{1}{c|}{$k$} & \multicolumn{1}{c|}{$N_k^*=(2k)!Q_k^*$} &
\multicolumn{1}{c|}{$Q_k^*=N_k^*/(2k)!$} &
\multicolumn{1}{c}{$\sqrt[k]{Q_k^*}$} \\ \hline 1 & 1 & $0.50000\dots$
& $0.50000\dots$ \\ 2 & 5 & $0.20833\dots$ & $0.45643\dots$ \\ 3 & 55
& $0.07638\dots$ & $0.42430\dots$ \\ 4 & 1023 & $0.02537\dots$ &
$0.39910\dots$ \\ 5 & 28207 & $0.00777\dots$ & $0.37854\dots$ \\ 6 &
1065317 & $0.00222\dots$ & $0.36129\dots$ \\
\end{tabular}
\end{center}
\bigskip

\noindent Using (\ref{wUBPn*}) for $k=6$, we get for each $\epsilon >0$

\[
P_n^*=\( \( \sqrt[6]{Q_6^*} +\epsilon \)^n \)=\( \( 0.361... +\epsilon \)^n \). \epf
\]

\bigskip

\begin{center}
{\Large {\bf Proof of Theorem \ref{ThmWBr}, lower bound.}}
\end{center}

To bound $P\( \pi\preceq \sigma \)$ from below we will use the criterion

\[
\pi\preceq \sigma\Longleftrightarrow E_i \(\pi \) \supseteq E_i \( \sigma\), \quad \forall \, i \le n.
\]

\noindent  First of all,

\bigskip

\begin{Lem} \label{LemEiContB} Let $i \in \left[n\right]$, $B \subseteq \left[ i-1 \right]$ ( $\left[0\right]=\emptyset$).  If $\pi \in S_n$ is chosen uniformly at random,  then

\[
P\( E_i \(\pi\) \supseteq B \) =  \frac{1}{|B|+1}. 
\]

\end{Lem}

\bpf By the definition of $E_i(\pi)$, 

\[
\{E_i(\pi)\supseteq B\}=\{\pi^{-1}(j)<\pi^{-1}(i),\,\forall\,j\in B\}.
\]

\noindent It remains to observe that $\pi^{-1}$ is also uniformly random. \epf

\bigskip

\noindent Lemma \ref{LemEiContB} implies the following key statement:

\bigskip

\begin{Lem} \label{LemEiContEi} Let $\pi,\sigma \in S_n$ be selected independently and uniformly at random. Then, for $i\in [n]$,

\[
P\( E_i \(\pi\) \supseteq E_i \(\sigma\) \) = H\(i\)/i,\quad H\(i\):=\sum_{j=1}^i\frac{1}{j}.
\]

\end{Lem}

\bpf By Lemma \ref{LemEiContB},

\begin{equation*} \begin{split}
P\(E_i\(\pi\)\supseteq E_i\(\sigma\)\)&=\sum_{B\subseteq
\left[i-1\right]}P\(E_i\(\pi\)\supseteq B\) P\(E_i\(\sigma\)=B\) \\
&=\sum_{B\subseteq \left[i-1\right]}\frac{P\(E_i\(\sigma\)=B\)}{|B|+1}
\\ &=E\left[ \, \frac{1}{|E_i(\sigma)|+1 } \, \right] \\
&=\sum_{j=0}^{i-1}\frac{1}{i \(j+1\)} =\frac{H\(i\)}{i}. \epf
\end{split} \end{equation*}

\bigskip

\noindent {\sc Note.~} In the second to last equality, we have used the 
fact that $|E_i(\sigma)|$ is distributed uniformly on
$\left\{0,1,\dots,i-1\right\}$. In addition, $|E_1(\sigma)|$,$\dots$, $|E_n(\sigma)|$ are
independent, a property we will use later. For completeness, here is
a bijective proof of these facts. By induction, the numbers $|E_i(\sigma)|$, $i\le t$,
determine uniquely the relative ordering of elements $1,\dots,t$ in the permutation $\sigma$.
Hence the numbers $|E_i(\sigma)|$, $i\in [n]$, determine $\sigma$ uniquely. Since the range of
$|E_i(\sigma)|$ is the set $\{0,\dots,i-1\}$ of cardinality $i$, and $|S_n|=n!$, it follows
that the numbers $|E_i(\sigma)|$, $i\in [n]$, are uniformly distributed, and independent of
each other.

\bigskip

Needless to say we are interested in $P(\pi\preceq \sigma)=P\left(\cap_{i=1}^n\{E_i(\pi)
\supseteq
E_i(\sigma)\}\right)$. Fortunately, the events $\{E_i(\pi)\supseteq E_i(\sigma)\}$ turn out to be
positively correlated, and  the product of the marginals $P(E_i(\pi)\supseteq
E_i(\sigma))$  bounds that probability from below.

\bigskip

\begin{Thm} \label{ThmWBrRest} Let $\pi, \sigma \in S_n$ be selected
 independently and uniformly at random.  Then
  
\begin{equation*} 
P \( \pi\preceq\sigma\)
\ge \prod_{i=1}^n P \( E_i \( \pi \) \supseteq E_i \(\sigma\) \) 
= \prod_{i=1}^n \frac{H\(i\)}{i}.
 \end{equation*}

\end{Thm}
 
\bpf  First notice that,
 conditioning on $\sigma$ and using the independence of $\pi$ and
 $\sigma$,
 
 \begin{equation*} \begin{split}
P \big( E_i \( \pi \) \supseteq & E_i \(\sigma\), \, \forall \, i \le
n \big) \\ &=E \left[ \, P \( E_i \( \pi \) \supseteq E_i \(\sigma\),
\, \forall \, i \le n \, | \, \sigma \, \)\right]=E \left[ \, P \( E_i
\( \pi \) \supseteq B_i, \, \forall \, i \le n \)|_{B_i=E_i
\(\sigma\)} \right].
  \end{split} \end{equation*}
 
 \noindent So our task is to bound $P \( E_i \( \pi \) \supseteq B_i,
 \, \forall \, i \le n \)$, where these $B_i$'s inherit the following
 property from the $E_i\(\sigma\)$'s:
 
 \[
 i \in E_j \( \sigma \) \textrm{ and } j \in E_k \( \sigma \)
 \Longrightarrow i \in E_k \(\sigma\).
 \]
 
 \bigskip
 
 \begin{Lem} \label{LemEiContBi} Let $n \ge 1$ be an integer, and let 
$B_i \subseteq \left[ n \right]$, $i=1,\dots,n$, be such that
 
 \[
 i\notin B_i\textrm{ and } i\in B_j,\, j \in B_k \Longrightarrow i \in B_k,\quad\forall\,i,j,k\in [n].
 \]
 
 \noindent Then, for $\pi \in S_n$ selected uniformly at random,  
 
 \[
 P \( E_i \( \pi \) \supseteq B_i, \, \forall \, i \le n \) \ge
 \prod_{i=1}^n \frac{1}{|B_i|+1}. 
 \]
 
 \end{Lem}

 \bigskip
 
 \noindent {\sc Proof of Lemma \ref{LemEiContBi}.~} Notice upfront that
 $\cup_iB_i\neq [n]$. Otherwise there would exist $i_1,\dots, i_s$ such
that $i_t\in B_{i_{t+1}}$, $1\le t\le s$, $(i_{s+1}=i_1)$, and -- using repeatedly the
property of the sets $B_i$ -- we would get that, say, $i_1\in B_{i_2}$
and $i_2\in B_{i_1}$, hence $i_2\in B_{i_2}$; contradiction.

Let $U_1,\dots,U_n$ be independent uniform-$[ 0,1 ]$ random
 variables. Let a random permutation $\omega$ be defined by

 \[
  \omega \(i \)=k 
 \Longleftrightarrow U_i \textrm{ is } k^\mathrm{th} \textrm{ smallest
 amongst } U_1,\dots,U_n.
 \]
 
 \noindent Clearly $\omega$ is distributed uniformly, and then so is $\pi:=\omega^{-1}$.
 With $\pi$ so defined, we obtain
 
 \[
 \left\{ E_i \( \pi \) \supseteq B_i, \, \forall \, i \le n \right\} =
 \left\{ \pi^{-1} \(i \) > \pi^{-1} \(j \), \, \forall \, j \in B_i,
 \, i \le n \right\} = \left\{ U_i > U_j, \, \forall \, j \in B_i,
 \, i \le n\right\}.
 \]
 
 \noindent Hence, the probability in question equals
 
 \[
 P \( U_i > U_j, \, \forall \, j \in B_i, \, i \le n \).
 \]
 
 \noindent We write this probability as the $n$-dimensional integral
 
 \begin{eqnarray*}
 &P \( U_i > U_j, \, \forall \, j \in B_i, \, i \le n\) =
 \idotsint\limits_D dx_1 \cdots dx_n,\\
 &D = \left\{ \(x_1,\dots,x_n \) \in \left[ 0,1 \right]^n \, : \,
 x_i>x_j, \, \forall \, j \in B_i, \, i \le n\right\}.
 \end{eqnarray*}
 
 \noindent Since $\cup_iB_i\neq  [n]$, we can choose an index $k \in \left[ n \right]$ such that $k
 \notin B_i$ for all $i$. Then we may rewrite the integral above as
 
\begin{eqnarray*}
&\int_0^1 \(\idotsint\limits_{D\(x_k\)} dx_1 \cdots dx_{k-1} dx_{k+1} \cdots
dx_n \) dx_k,\\
&D\(x_k\)= \left\{ \(x_1,\dots,x_{k-1},x_{k+1},\dots,x_n \) \in \left[
 0,1 \right]^{n-1} \, : \, x_i>x_j, \, \forall \, j \in B_i, \, i \le n\right\}.
\end{eqnarray*}

\noindent On $D\(x_k\)$, the only inequalities involving $x_k$ are
of the form $x_k > x_j$, $j\in B_k$. This suggests scaling those $x_j$ by $x_k$,
i.e. introducing new variables
$t_j:=x_j/x_k$, so that $t_j\in \left[0,1\right]$, $j\in B_k$.  To keep
notation uniform, let us also replace the remaining $x_i$, $i\notin B_k\cup\{k\}$, 
with $t_i$.  Let
$\mathfrak{D}\(x_k\)$ denote the integration region for the new
variables $t_i$, $i\neq k$. Explicitly, the constraints 
$x_j < x_k$, $j\in B_k$, become $t_j< 1$, $j\in B_k$. Obviously each listed
constraint $x_a < x_b$ ($a,b\in B_k$) is replaced, upon scaling, with $t_a<t_b$.
We only rename the other variables, so every constraint $x_a<x_b$ 
($a,b\notin B_k$) similarly becomes $t_a<t_b$. By the property of the
sets $B_i$, there are no inequalities $x_a > x_b$, $a\in B_k$, $ b\notin B_k$ (since the 
presence of this inequality implies $b \in B_a$).
The only remaining inequalities are all of the type $x_a < x_b$, $a\in B_k$, 
$b\notin B_k$. In the new variables, such a constraint becomes $x_kt_a < t_b$, and
it is certainly satified if $t_a < t_b$, as $x_k\le 1$. Hence, $\mathfrak{D}\(x_k\)
 \supseteq D^*$, where

\[
D^*:=\left\{ \(t_1,\dots,t_{k-1},t_{k+1},\dots,t_n \) \in [0,1]^{n-1}
\, : \, t_i > t_j, \, \forall \, j \in B_i, \, i \neq k \right\},
\]

\noindent and $D^*$ does not depend on $x_k$! Observing that the constraints
that determine $D^*$ are those for $D$ with the constraints $x_i<x_k$, $i
\in B_k$, removed, we conclude that the innermost integral over $D(x_k)$ is
bounded below by $x_k^{|B_k|}P(U_i>U_j,\, \forall \, j \in B_i, \, i \neq k)$. 
($x_k^{|B_k|}$ is the Jacobian of the linear transformation
$\{x_i\}_{i\neq k}\to \{t_i\}_{i\neq k}$.) Integrating with respect to $x_k$,
we arrive at

 \begin{equation} \label{IndPUU>}
 P \( U_i > U_j, \, \forall \, j \in B_i, \, i \le n \)\ge 
 \frac{1}{|B_k|+1} \cdot P \( U_i > U_j, \, \forall \, j \in B_i, \, i \neq k\).
 \end{equation}
 
 \noindent By induction on the number of sets $B_i$, with Lemma \ref{LemEiContB} 
 providing basis of induction and (\ref{IndPUU>}) -- the inductive step, we get

 \begin{equation*}
 P \( U_i > U_j, \, \forall \, j \in B_i, \, i \le n\) \ge
 \prod_{i=1}^n \frac{1}{|B_i|+1}. \epf
 \end{equation*}
 
 \bigskip
 
 \noindent The rest is short. First, by Lemma \ref{LemEiContBi},

\begin{equation*} \begin{split}
P(E_i(\pi)\supseteq E_i(\sigma),\,\forall \, i\le n)&=E\left[\left.
P(E_i(\pi)\supseteq B_i,\,\forall \, i\le n)\right|_{B_i=E_i(\sigma)}\right]\\
&\ge E\left[  \prod_{i=1}^n \frac{1}{|E_i \( \sigma \)| +1} \right].
\end{split} \end{equation*}
 
\noindent Since the cardinalities  $|E_i \(\sigma \)|$ are independent, the last
expected value equals
 
 \begin{equation*}
 \prod_{i=1}^n E \left[ \, \frac{1}{|E_i \( \sigma \)| +1} \, \right]
 = \prod_{i=1}^n \(\frac{1}{i} \sum_{j=0}^{i-1}\frac{1}{j+1}\) =
 \prod_{i=1}^n \frac{H\(i\)}{i};
 \end{equation*}
 
\noindent for the second to last equality see the proof of Lemma \ref{LemEiContEi}. \epf
 
 \bigskip

 \noindent {\sc Note.~} Let ${\cal P}$ be a
 poset on $\left[ n \right]$, and put $B_i:= \big\{ j \in {\cal P} \,
 : \, j < i \text{ in } {\cal P}\big\}$. $B_i \cup \left\{i\right\}$
 is called the \emph{order ideal at $i$}. By the properties of ${\cal
 P}$, the $B_i$'s satisfy the hypotheses of Lemma \ref{LemEiContBi},
 so letting $e\({\cal P}\)$ denote the number of linear extensions of
 ${\cal P}$ we get
 
 \begin{equation*} \begin{split}
 P \( E_i \( \pi \) \supseteq B_i, \, \forall \, i \le n
 \)&=\frac{|\left\{ \omega \, : \, \omega\(i\)>\omega\(j\),\, \forall \, j \in B_i, \, i \le n\right\}|}{n!}\\
&=\frac{e\({\cal
 P}\)}{n!}  \ge \prod_{i=1}^n \frac{1}{|B_i|+1}.
 \end{split} \end{equation*}
 
 \noindent Thus we have proved
 
\bigskip 
 
\begin{Cor} \label{LinExtCor} For a poset ${\cal P}$ with $n$ elements,

\[
e\({\cal P}\) \ge n!\Big/\prod\limits_{i=1}^n d\(i\),\quad 
d\(i\):=|\left\{ j \in {\cal P} \, : \, j \le i
\text{ in } {\cal P}\right\}|. \epf
\]
 
\end{Cor}

\bigskip

\noindent In a very special case of ${\cal P}$, whose Hasse diagram is a forest of rooted
trees with edges directed away from the roots, this simple bound is actually the
value of $e({\cal P})$ (see \cite[Section 5.1.4, ex. 20]{K}). There exist bounds
for the number of linear extensions in the case of the Boolean lattice (see Brightwell 
and Tetali \cite{BT}, Kleitman and Sha \cite{KS}), but the lower bound in
Corollary \ref{LinExtCor} seems to be new.

\bigskip

\begin{center}
{\Large {\bf Numerics}}
\end{center}

From computer-generated data we have collected, it
appears that our $O\(n^{-2}\)$ upper bound correctly predicts the qualitative behavior of $P\( \pi \leq \sigma \)$. The data suggests that $P\( \pi \leq
\sigma \)$ is of exact order $n^{-\(2+\delta \)}$ for some $\delta \in
\left[0.5,1\right]$, which begs the question of how to improve on our current bound. Writing 
$P_n=P\( \pi \leq \sigma \)$, below is a graph
(based on this numerical experimentation) exhibiting convergence to
the exponent $-a$ in the asymptotic equation $P_n\sim cn^{-a}$, $c>0$ a
constant, and $-a$ appears to be near $-2.5$:

\bigskip
\begin{center}
\includegraphics[scale=0.7]{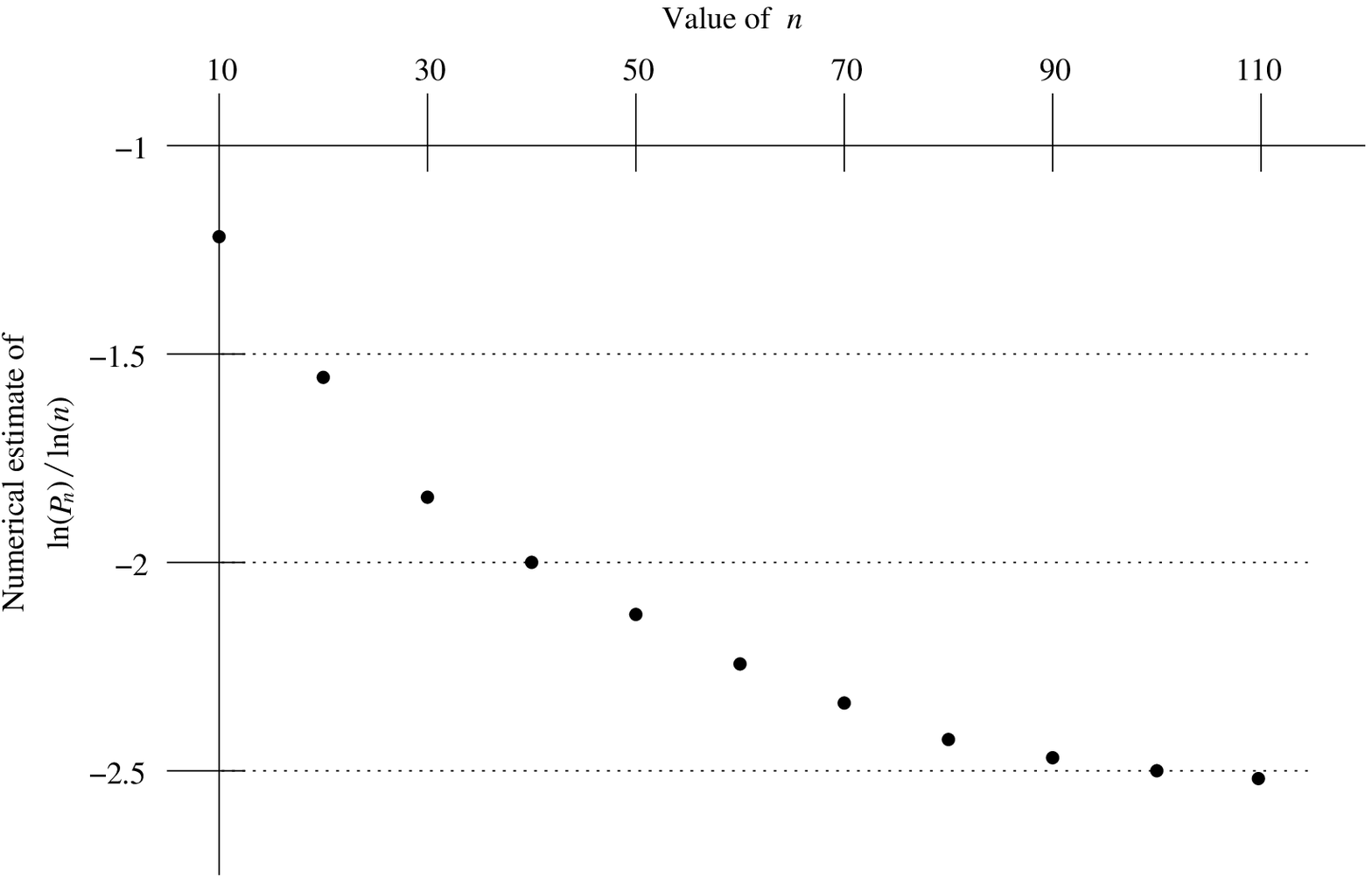}
\end{center}
\bigskip

\noindent Here is a portion of the accompanying table used to generate
this graph:

\bigskip
\begin{center}
\begin{tabular}{r|r|l|l}
\multicolumn{1}{c|}{$n$} & \multicolumn{1}{c|}{$R_n$} &
\multicolumn{1}{c|}{Estimate of $P_n \approx \frac{R_n}{10^9}$} &
\multicolumn{1}{c}{Estimate of $\ln (P_n)/\ln n$} \\ \hline 10 &
61589126 & $0.0615891\dots$ & $-1.21049\dots$ \\ 30 & 1892634 &
$0.0018926\dots$ & $-1.84340\dots$ \\ 50 & 233915 & $0.0002339\dots$ &
$-2.13714\dots$ \\ 70 & 50468 & $0.0000504\dots$ & $-2.32886\dots$ \\
90 & 14686 & $0.0000146\dots$ & $-2.47313\dots$ \\ 110 & 5174 &
$0.0000051\dots$ & $-2.58949\dots$ \\
\end{tabular}
\end{center}
\bigskip

\noindent Here, $R_n$ is the number of pairs $(\pi,\sigma)$ out of
$10^9$ randomly-generated pairs such that we had $\pi \le \sigma$. We
have also utilized the computer to find the actual probability $P_n$
for $n=1,2,\dots,9$.  Below is a table of these true proportions:

\bigskip
\begin{center}
\begin{tabular}{r|r|l}
\multicolumn{1}{c|}{$n$} & \multicolumn{1}{c|}{$(n!)^2P_n$} & \multicolumn{1}{c}{$P_n$} \\ 
\hline 1 & $1$ & $1.00000\dots$ \\ 
2 & $3$ & $0.75000\dots$ \\ 
3 & $19$ & $0.52777\dots$ \\ 
4 & $213$ & $0.36979\dots$ \\ 
5 & $3781$ & $0.26256\dots$ \\ 
6 & $98407$ & $0.18982\dots$ \\ 
7 & $3550919$ & $0.13979\dots$ \\ 
8 & $170288585$ & $0.10474\dots$ \\ 
9 & $10501351657$ & $0.07974\dots$ \\
\end{tabular}
\end{center}
\bigskip

Concerning the weak Bruhat order, computer-generated data suggests
that $P\( \pi \preceq \sigma \)$ is of exact order $(0.3)^n$. So our current 
upper bound $O((0.362)^n)$ is a qualitative match for $P\( \pi \preceq \sigma \)$, 
but it appears that improvements are possible here also.  Writing $P_n^*=P\(
\pi \preceq \sigma \)$, below is a graph (based on our numerical
experiments) exhibiting convergence to the ratio $\rho$ in the
asymptotic equation $P_n^* \sim c\rho^n$, $c>0$ a constant, and $\rho$
appears to be near $0.3$:

\bigskip
\begin{center}
\includegraphics[scale=.7]{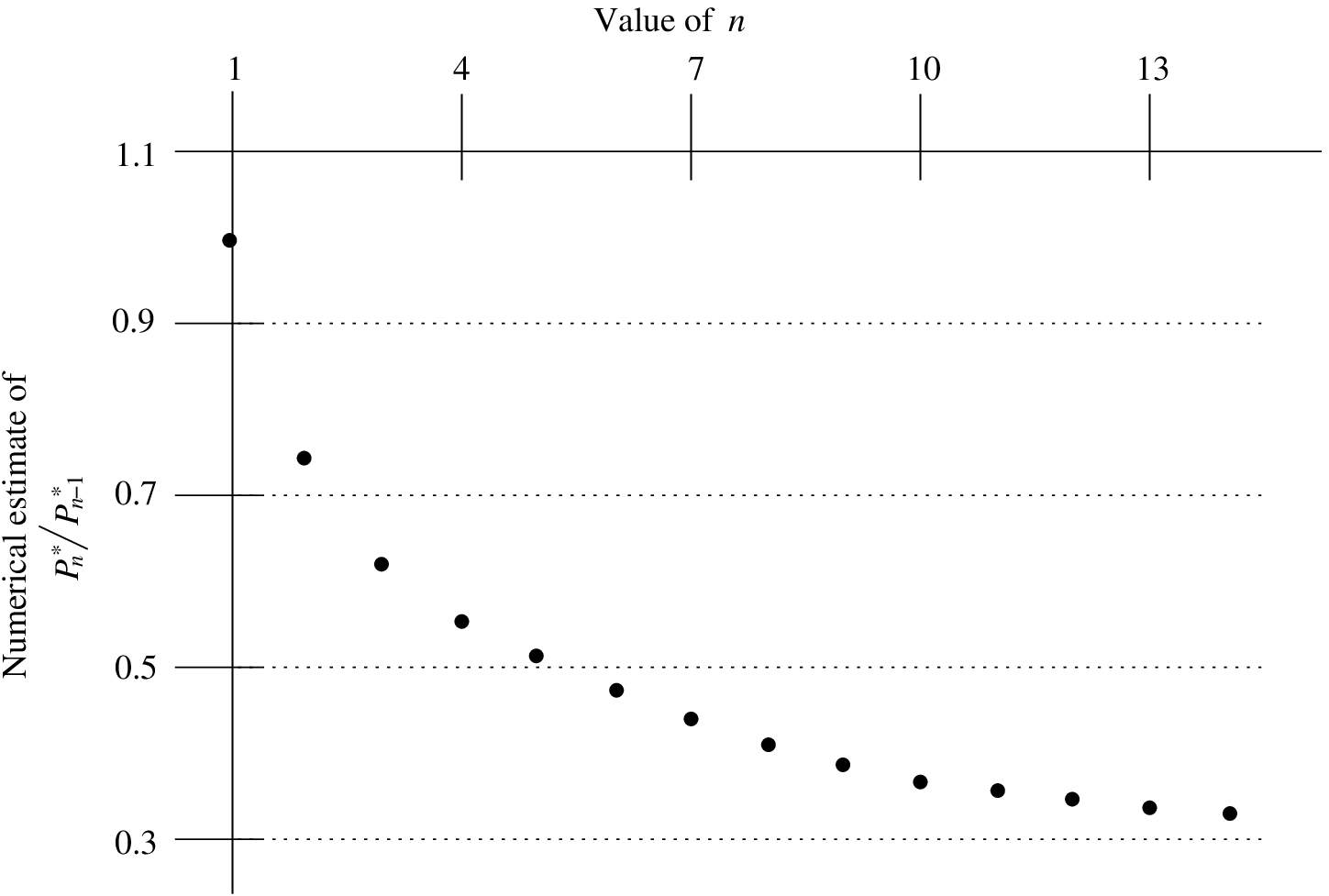}
\end{center}
\bigskip

\noindent Here is a portion of the accompanying table used to generate
this graph:

\bigskip
\begin{center}
\begin{tabular}{r|r|l|l}
\multicolumn{1}{c|}{$n$} & \multicolumn{1}{c|}{$R_n^*$} &
\multicolumn{1}{c|}{Estimate of $P_n^* \approx \frac{R_n^*}{10^9}$} &
\multicolumn{1}{c}{Estimate of $P_n^* /P_{n-1}^*$} \\ \hline 10 &
1538639 & $0.0015386\dots$ & $0.368718\dots$ \\ 11 & 541488 &
$0.0005414\dots$ & $0.351926\dots$ \\ 12 & 184273 & $0.0001842\dots$ &
$0.340308\dots$ \\ 13 & 59917 & $0.0000599\dots$ & $0.325153\dots$ \\
14 & 18721 & $0.0000187\dots$ & $0.312448\dots$ \\ 15 & 5714 &
$0.0000057\dots$ & $0.305218\dots$ \\ 16 & 1724 & $0.0000017\dots$ &
$0.301715\dots$ \\
\end{tabular}
\end{center}
\bigskip

\noindent Here, $R_n^*$ is defined analogously to $R_n$ above. Below
is a table with the true proportion $P_n^*$ for $n=1,2,\dots,9$:

\bigskip
\begin{center}
\begin{tabular}{r|r|l}
\multicolumn{1}{c|}{$n$} & \multicolumn{1}{c|}{$(n!)^2P_n^*$} & \multicolumn{1}{c}{$P_n^*$} \\ 
\hline 1 & $1$ & 1.00000\dots \\ 
2 & $3$ & 0.75000\dots \\ 
3 & $17$ & 0.47222\dots \\ 
4 & $151$ & 0.26215\dots \\ 
5 & $1899$ & 0.13187\dots \\ 
6 & $31711$ & 0.06117\dots \\ 
7 & $672697$ & 0.02648\dots \\ 
8 & $17551323$ & 0.01079\dots \\ 
9 & $549500451$ & 0.00417\dots \\
\end{tabular}
\end{center}
\bigskip

\noindent Surprisingly, our Theorem \ref{ThmWBr} lower bound for
$P_n^*$ is quite good for these smallish values of $n$:

\bigskip
\begin{center}
\begin{tabular}{r|r|l}
\multicolumn{1}{c|}{$n$} & \multicolumn{1}{c|}{$(n!)^2\prod_{i=1}^n \(H\(i\)/i\)$} & \multicolumn{1}{c}{$\prod_{i=1}^n \(H\(i\)/i\)$} \\ 
\hline 1 & $1.0\dots$ & $1.00000\dots$ \\ 
2 & $3.0\dots$ & $0.75000\dots$ \\ 
3 & $16.5\dots$ & $0.45833\dots$ \\ 
4 & $137.5\dots$ & $0.23871\dots$ \\ 
5 & $1569.8\dots$ & $0.10901\dots$ \\ 
6 & $23075.9\dots$ & $0.04451\dots$ \\ 
7 & $418828.3\dots$ & $0.01648\dots$ \\ 
8 & $9106523.1\dots$ & $0.00560\dots$ \\ 
9 & $231858583.9\dots$ & $0.00176\dots$ \\
\end{tabular}
\end{center}
\bigskip

\begin{center}
{\Large {\bf Acknowledgments}}
\end{center}

This work was inspired by a thought-provoking talk Mark Skandera gave
at the MIT Combinatorics Conference honoring Richard Stanley's
60$^\mathrm{th}$ birthday. We are grateful to Mark for an enlightening
follow-up discussion of comparability criteria for the Bruhat
order. We thank Sergey Fomin for encouragement and for introducing us
to an instrumental notion of the permutation-induced poset. Without Ed Overman's 
invaluable guidance we would not have been able to 
obtain the numerical results in this paper. Craig Lennon gave us an idea for proving 
an exponential lower bound in the
case of strong Bruhat order.

\pagebreak

\begin{center}

\end{center}

\vfill
\end{document}